\documentclass[final,hidelinks,onefignum,onetabnum]{siamart251216}

\usepackage[utf8]{inputenc}
\usepackage{amssymb,amsfonts,mathrsfs}
\usepackage{booktabs,longtable,tabularx,siunitx,multirow}
\usepackage{graphicx,float,tikz}

\newsiamremark{remark}{Remark}

\headers{A Geometry-Robust Ne\v{c}as-Type Inequality}{Q. Xin, Y. Xie, and S. Gong}
\title{A Geometry-Robust Ne\v{c}as-Type Inequality for Stokes Preconditioning}
\author{Qi Xin\thanks{School of Science and Engineering, The Chinese
University of Hong Kong, Shenzhen, Guangdong 518172, China;
Shenzhen International Center for Industrial and Applied Mathematics,
Shenzhen Research Institute of Big Data, Shenzhen, Guangdong 518172, China
(\email{qixin1@link.cuhk.edu.cn}).}
\and Yan Xie\thanks{Center for Engineering Computing and Optimization,
Shenzhen Research Institute of Big Data, Shenzhen, Guangdong
518172, China (\email{xieyan@sribd.cn}).}
\and Shihua Gong\thanks{Corresponding author. School of Science and Engineering,
The Chinese University of Hong Kong, Shenzhen, Guangdong 518172, China;
Shenzhen International Center for Industrial and Applied Mathematics,
Shenzhen Research Institute of Big Data, Shenzhen, Guangdong 518172, China;
Shenzhen Loop Area Institute, Shenzhen, Guangdong 518040, China
(\email{gongshihua@cuhk.edu.cn}).}}

\newcommand{\R}{\mathbb{R}}
\newcommand{\eps}{\epsilon}

\DeclareMathOperator{\dist}{dist}

\begin{document}

\maketitle

\begin{abstract}
Pressure mass preconditioners for the incompressible Stokes equations can deteriorate on domains with narrow gaps. We develop a pressure preconditioner that uses the distance to the boundary to account for spatially varying geometric scales. Its application consists of a pressure mass solve and a weighted pressure Poisson solve, with the squared boundary distance as the diffusion coefficient. The theoretical basis is a distance-weighted pressure norm defined by splitting the pressure into an $L^2$ component and a component controlled by its weighted gradient. We prove that this norm is bounded by the negative norm of the pressure gradient on every bounded connected open set, with a constant depending only on the dimension. A reverse bound follows from the boundary Hardy inequality, giving norm equivalence uniformly over families with bounded Hardy constants. Two-dimensional finite element experiments on slender, perforated, and composite domains demonstrate stable performance across the tested geometric parameters and mesh refinements. The smallest eigenvalues of the preconditioned pressure Schur complement remain close to $0.54$, and block preconditioning with algebraic multigrid auxiliary solves requires about 20--30 FGMRES iterations across the geometric sweeps.
\end{abstract}



\section{Introduction}\label{sec:intro}

We study the efficient numerical solution of the steady Stokes equations
on domains with challenging geometries. For unit viscosity and homogeneous
Dirichlet velocity data, the equations take the form
\begin{equation}\label{eq:intro_stokes}
    -\Delta\mathbf u+\nabla p=\mathbf f,\qquad
    \nabla\cdot\mathbf u=0\quad\text{in }\Omega,\qquad
    \mathbf u=\mathbf0\quad\text{on }\partial\Omega,
\end{equation}
where $\Omega\subset\R^d$, $d\ge2$, is a bounded connected Lipschitz
domain, $\mathbf u$ is the velocity, and $p$ is the pressure, normalized
to have zero mean. Mixed finite-element discretizations lead to
saddle-point systems, for which effective pressure preconditioning is
essential to efficient iterative solution \cite{Elman2014,Benzi2005}.

The viscous term provides direct control of
$\|\nabla\mathbf u\|_{L^2(\Omega)}$, while the pressure acts as a Lagrange
multiplier enforcing incompressibility and has no coercive term of its
own. Pressure stability therefore depends on its coupling to
the velocity through the divergence operator. The corresponding inf-sup
condition requires that the $L^2$ norm of a mean-zero pressure be
controlled by the negative norm of its gradient. Equivalently, one needs
the classical Ne\v{c}as inequality \cite{Necas1967,Girault1986,Bramble2003}:
\begin{equation}\label{eq:intro_classical_necas}
    \|q\|_{L^2(\Omega)}
    \le C_N(\Omega)\,\|\nabla q\|_{-1,\Omega},
    \qquad q\in L_0^2(\Omega),
\end{equation}
where $L_0^2(\Omega)$ is the space of square-integrable functions with
zero mean, and the negative norm is dual to the gradient norm on
$H_0^1(\Omega)^d$. This estimate also underpins the comparison of the
pressure Schur complement with a pressure mass operator in block
preconditioners \cite{Elman2014,loghin2004analysis}.

In the classical negative-norm theory, the Ne\v{c}as inequality is also
stated with a lower-order term,
\begin{equation*}
    \|q\|_{L^2(\Omega)}
    \le C(\Omega)\bigl(\|q\|_{H^{-1}(\Omega)}
                       +\|\nabla q\|_{-1,\Omega}\bigr),
    \qquad q\in L^2(\Omega).
\end{equation*}
On a bounded connected Lipschitz domain, a compactness argument removes
the lower-order term after fixing the mean, giving
\eqref{eq:intro_classical_necas}; see \cite{Necas1967,Girault1986}.
Bramble \cite{Bramble2003} gives a direct proof on Lipschitz domains
and explains its equivalence to the Stokes inf-sup condition.
Amrouche, Ciarlet and Mardare \cite{AmroucheCiarletMardare2015}
place these estimates in a broader equivalence theorem linking the
Lions lemma, the closed range of the gradient, solvability of the
divergence equation, and a weak form of de Rham's theorem.

The classical Ne\v{c}as constant depends on the domain geometry and may
become unbounded along degenerating families. This dependence has been
studied extensively through the equivalent
Ladyzhenskaya--Babu\v{s}ka--Brezzi (LBB) formulation.
Chizhonkov and Olshanskii \cite{Chizhonkov2000} obtained explicit,
asymptotically accurate estimates on model domains such as strips and
rings, and related the continuous constant to discrete stability.
Dobrowolski \cite{Dobrowolski2003} established aspect-ratio estimates on
stretched domains; his subsequent work \cite{Dobrowolski2012} used coarse
pressure spaces to construct preconditioners with bounds independent of
the aspect ratio. Our preprint \cite{infsup_perforation} gives a sharp
LBB estimate for tight domains with periodic pillars.
These results quantify how geometric degeneration enters
the classical pressure estimate and how its effect can be compensated.

To obtain uniform stability estimates, one can incorporate the geometric scale into the pressure norm itself. For slender domains of width $\eps$, Sande, Koch,
Kuchta and Mardal \cite{sande2025robust} established a uniform equivalence
between $\|\nabla q\|_{-1,\Omega}$ and the $K$-functional norm
$\|q\|_{L^2+\eps H^1}$, using a coarse partition and quasi-interpolation.
For periodic pillar arrays with hole size proportional to the period $\eps$,
the sharp LBB estimate of \cite{infsup_perforation} is complemented by
the norm equivalence
$\|\nabla q\|_{-1,\Omega}\asymp\|q\|_{L^2+\sigma_\eps H^1}$ in
\cite{xin2026uniform}, where $\sigma_\eps\asymp\eps$ in this regime.
The latter proof combines a weighted Ne\v{c}as estimate on a fixed filled
domain with a dual pressure extension built from the restriction operators
of \cite{Allaire1991a,Allaire1991,Lu2020}; its constants are uniform in
the period for a fixed macroscopic domain and reference geometry.
In Ne\v{c}as language, these norm characterizations provide pressure
estimates adapted to the scale of the specified geometric family.

Spatially varying geometric scales are common in practical flow domains,
including branching vascular geometries \cite{Figueroa2006},
microfluidic devices employing deterministic lateral displacement (DLD)
arrays \cite{Loutherback2012,Shen2026PredictiveMicrofluidic}, and the
pore spaces of natural rocks \cite{Mostaghimi2013}. Sande et al.\
\cite[Remark~3 and Example~5]{sande2025robust} also discuss a variable
local width and the associated weighted pressure operator. The analytical
question pursued here is to identify a canonical local scale and prove
the corresponding Ne\v{c}as-type inequality on arbitrary bounded connected
open sets, with a constant depending only on the dimension. The boundary
distance supplies such a scale intrinsically, without a prescribed channel
width, period, or coarse partition.

Let $d\ge2$, let $\Omega\subset\R^d$ be a bounded connected open set, and set
$\delta(x)=\dist(x,\partial\Omega)$. Let $\mathcal H_\delta(\Omega)$
consist of the functions $\phi\in L^2(\Omega)\cap H^1_{\mathrm{loc}}(\Omega)$
with $\delta\nabla\phi\in L^2(\Omega)^d$. Define the distance-weighted
$K$-functional norm by
\begin{equation}
\label{eq:intro_K_delta}
    \|q\|_{K_\delta}
    \;:=\;
    \inf_{\phi\in\mathcal H_\delta(\Omega)}
    \bigl(
        \|q-\phi\|_{L^2(\Omega)}^2
        +\|\delta\nabla\phi\|_{L^2(\Omega)}^2
    \bigr)^{1/2}.
\end{equation}
Our main result, Theorem~\ref{thm:weighted_necas}, is the inequality
\begin{equation}\label{eq:intro_weighted_necas}
    \|q\|_{K_\delta}
    \le C_d\,\|\nabla q\|_{-1,\Omega},
    \qquad q\in L_0^2(\Omega),
\end{equation}
with $C_d$ depending only on the dimension. The geometric information is
carried by $\delta$ in the norm, while the constant is independent of
aspect ratios, perforations, and contrasts between local length scales.
Since $\|q\|_{K_\delta}\le\|q\|_{L^2(\Omega)}$, it controls pressure in
a weaker, geometrically adapted norm.

The proof starts from the classical Ne\v{c}as inequality on a cube,
which controls a function up to its mean with a constant independent of
the cube size. A Whitney decomposition allows us to apply this estimate
on slightly enlarged cubes whose sizes are comparable to the local
boundary distance. Each cube has its own mean, and these local constants
need not agree. We blend them into a smooth function $\phi$ using a
partition of unity and use $\phi$ as a competitor in
\eqref{eq:intro_K_delta}. Both the error $q-\phi$ and the weighted
gradient $\delta\nabla\phi$ are controlled by the local oscillations:
the factor $\delta$ compensates for the derivatives of the partition
functions. Bounded overlap and negative-norm localization then give
\eqref{eq:intro_weighted_necas}.

If $\Omega$ also satisfies the boundary Hardy inequality
$\|v/\delta\|_{L^2(\Omega)}\le C_H\|\nabla v\|_{L^2(\Omega)}$ for
$v\in H_0^1(\Omega)$, then integration by parts gives a reverse
estimate. Combining this estimate with \eqref{eq:intro_weighted_necas}
yields the norm
equivalence in Corollary~\ref{cor:delta_equivalence}:
\begin{equation}\label{eq:intro_delta_equivalence}
    C_d^{-1}\,\|q\|_{K_\delta}
    \le \|\nabla q\|_{-1,\Omega}
    \le \sqrt{1+C_H^2}\,\|q\|_{K_\delta},
    \qquad q\in L_0^2(\Omega).
\end{equation}
The equivalence is uniform on any family with uniformly bounded Hardy
constants. For example, Ancona's bound $C_H\le4$ for simply connected
planar domains \cite{Ancona1986,LaptevSobolev2008} gives the reverse
constant $\sqrt{17}$, independently of channel widths and their contrast;
see Remark~\ref{rem:uniform_hardy}. Thus the additional geometric
dependence in the two-sided norm characterization is isolated in the
Hardy constant.

For the Stokes problem, this norm equivalence provides an analytical
basis for pressure preconditioning through the operator-preconditioning framework of
Mardal--Winther \cite{Mardal2010}, building on Bramble--Pasciak
\cite{Bramble1988} and Loghin--Wathen \cite{loghin2004analysis}.
The squared negative norm of the pressure gradient is the quadratic form
of the continuous pressure Schur complement, so an equivalent pressure
norm supplies a corresponding Riesz map. The distance-weighted
$K$-functional motivates an inverse
pressure block combining a mass inverse and a
boundary-distance-weighted stiffness inverse, with no prescribed global
length scale. Transferring uniform continuous bounds to a discretization
requires suitable uniform discrete stability and comparison of the
discrete weighted norm with $K_\delta$.
Numerical results on slender, perforated, and composite domains show
that the computed smallest eigenvalues of the preconditioned pressure
Schur complement remain uniformly bounded away from zero across the
tested geometric variations, consistent with the geometry-independent
lower bound in the continuous estimate.

The remainder of the paper is organized as follows.
Section~\ref{sec:problem} introduces the negative norm of the gradient
and the distance-weighted $K$-functional norm.
Section~\ref{sec:main} proves the geometry-robust weighted Ne\v{c}as
inequality through local estimates on Whitney cubes and a smooth
partition of unity. Section~\ref{subsec:norm_equivalence} then uses the
boundary Hardy inequality to derive the norm equivalence as a
corollary. Section~\ref{sec:preconditioning} develops the Stokes
preconditioning application, and Section~\ref{sec:numerical} reports
numerical experiments.
Section~\ref{sec:conclusion} concludes.

\section{Mathematical Setting}\label{sec:problem}

Let $\Omega\subset\R^d$, $d\ge2$, be a bounded connected open set.
Write $L_0^2(\Omega)=\{q\in L^2(\Omega):\int_\Omega q\,dx=0\}$.
For $q\in L^2(\Omega)$, define the negative norm of its gradient by
\begin{equation}\label{eq:gradient_negative_norm}
    \|\nabla q\|_{-1,\Omega}
    :=\sup_{\mathbf v\in H_0^1(\Omega)^d\setminus\{\mathbf0\}}
      \frac{|\langle\nabla q,\mathbf v\rangle|}
           {\|\nabla\mathbf v\|_{L^2(\Omega)}},
\end{equation}
where $\langle\nabla q,\mathbf v\rangle=-\int_\Omega q\,\nabla\cdot\mathbf v\,dx$.

To define the geometrically adapted norm, let
\begin{equation}\label{eq:delta_def}
    \delta(x):=\dist(x,\partial\Omega),
    \qquad x\in\Omega.
\end{equation}
The natural weighted energy space is
\begin{equation}\label{eq:weighted_energy_space}
    \mathcal H_\delta(\Omega)
    :=
    \left\{
      \phi\in L^2(\Omega):
      \phi\in H^1_{\mathrm{loc}}(\Omega),\
      \delta\nabla\phi\in L^2(\Omega)^d
    \right\}.
\end{equation}
We define the distance-weighted $K$-functional norm by
\begin{equation}\label{eq:K_delta_def}
    \|q\|_{K_\delta}
    :=
    \inf_{\phi\in\mathcal H_\delta(\Omega)}
    \Bigl(
        \|q-\phi\|_{L^2(\Omega)}^2
        +\|\delta\nabla\phi\|_{L^2(\Omega)}^2
    \Bigr)^{1/2},
    \qquad q\in L^2(\Omega).
\end{equation}

The functional is unchanged when a constant is added to $q$: the change
of variables $\phi\mapsto\phi+c$ shows
$\|q+c\|_{K_\delta}=\|q\|_{K_\delta}$.  Thus it is naturally a norm on
$L^2(\Omega)/\mathbb R$.  For completeness, if
$\|q\|_{K_\delta}=0$, there are $\phi_n\in\mathcal H_\delta(\Omega)$ such that
$\phi_n\to q$ in $L^2$ and $\delta\nabla\phi_n\to0$ in $L^2$.  On every
$U\Subset\Omega$, $\delta$ is bounded below, and hence
$\nabla\phi_n\to0$ in $L^2(U)$.  It follows distributionally that
$\nabla q=0$ in $\Omega$; connectedness implies that $q$ is constant.
Therefore \eqref{eq:K_delta_def} is a genuine norm on $L_0^2(\Omega)$.

\section{A Geometry-Robust Weighted Ne\v{c}as Inequality}\label{sec:main}

We prove the following distance-weighted inequality.

\begin{theorem}[Geometry-robust distance-weighted Ne\v{c}as-type inequality]
\label{thm:weighted_necas}
Let $d\ge2$, let $\Omega\subset\R^d$ be a bounded connected open set,
and set $\delta(x)=\dist(x,\partial\Omega)$. Then
\begin{equation}\label{eq:weighted_necas}
    \|q\|_{K_\delta}
    \;\le\;
    C_d\,\|\nabla q\|_{-1,\Omega},
    \qquad q\in L_0^2(\Omega).
\end{equation}
The constant $C_d>0$ depends only on $d$, independently of the domain.
\end{theorem}

The proof combines local Ne\v{c}as estimates on Whitney cubes with a
smooth partition of unity. We first establish the covering and local
estimates, then construct a competitor for the infimum in
\eqref{eq:K_delta_def} to obtain \eqref{eq:weighted_necas}.
All constants in the proof depend only on $d$.

\subsection{Whitney covering and partition of unity}
\label{subsec:whitney_covering}

\begin{lemma}[Whitney covering and partition of unity]
\label{lem:whitney_covering}
There is a countable family $\mathcal W$ of closed cubes with disjoint
interiors whose union is $\Omega$.
Write $\ell_Q$ for the side length
of $Q\in\mathcal W$ and set $U_Q=\operatorname{int}(\tfrac98Q)$,
where the dilation is concentric.
There are positive constants $c_d,C_d,N_d,A_d,B_d$, depending only on
$d$, such that the following properties hold.
The sets $U_Q\Subset\Omega$ form a
locally finite open cover, and
\begin{equation}\label{eq:whitney_covering_bounds}
 c_d\ell_Q\le\delta(x)\le C_d\ell_Q
 \quad(x\in U_Q),\qquad
 \sum_{Q\in\mathcal W}\chi_{U_Q}\le N_d.
\end{equation}
There are functions $\psi_Q\in C_c^\infty(U_Q)$ such that
\begin{equation}\label{eq:whitney_partition}
 0\le\psi_Q\le1,\qquad
 \sum_{Q\in\mathcal W}\psi_Q=1\text{ in }\Omega,\qquad
 |\nabla\psi_Q|\le A_d\ell_Q^{-1}.
\end{equation}
In particular, $\delta|\nabla\psi_Q|\le B_d$ in $\Omega$.
\end{lemma}

\begin{proof}
The Whitney covering and the bounds in
\eqref{eq:whitney_covering_bounds}--\eqref{eq:whitney_partition}
are standard; see \cite[Chapter~VI, Section~1]{Stein1970}, in particular
Section~1.3 for the partition of unity.
Combining $\delta\le C_d\ell_Q$ on $U_Q$ with
$|\nabla\psi_Q|\le A_d\ell_Q^{-1}$ gives
$\delta|\nabla\psi_Q|\le B_d:=C_dA_d$ there; outside $U_Q$,
$\nabla\psi_Q=0$.
\end{proof}

\subsection{Local Ne\v{c}as estimates and negative-norm localization}
\label{subsec:local_estimates}

For an open set $U\subset\Omega$, let $\|\nabla q\|_{-1,U}$ denote
the negative norm in \eqref{eq:gradient_negative_norm} with $\Omega$
replaced by $U$.

\begin{lemma}[Local Ne\v{c}as estimate on a cube]\label{lem:cube_necas}
For every open cube $U\subset\R^d$ and every $r\in L^2(U)$,
\begin{equation}\label{eq:cube_necas}
 \|r-r_U\|_{L^2(U)}\le C_d\|\nabla r\|_{-1,U},
 \qquad r_U:=\frac1{|U|}\int_Ur\,dx.
\end{equation}
The constant is independent of the position and side length of $U$.
\end{lemma}

This is the classical Ne\v{c}as inequality \cite{Necas1967}, scaled
from a reference cube.

\begin{lemma}[Negative-norm localization]\label{lem:dual_localization}
Let $\mathcal U$ be a finite or countable family of open subsets of
$\Omega$ with $\sum_{U\in\mathcal U}\chi_U\le M$ almost everywhere.
Then, for every $q\in L^2(\Omega)$,
\begin{equation}\label{eq:dual_localization}
 \sum_{U\in\mathcal U}\|\nabla q\|_{-1,U}^2
 \le M\|\nabla q\|_{-1,\Omega}^2.
\end{equation}
\end{lemma}

\begin{proof}
Fix a finite subfamily and $0<\eta<1$, and write
$a_U=\|\nabla q\|_{-1,U}$. Omitting terms with $a_U=0$, choose
$\mathbf v_U\in H_0^1(U)^d$ with
$\|\nabla\mathbf v_U\|_{L^2(U)}=1$ and
$\int_Uq\,\nabla\cdot\mathbf v_U\,dx\ge(1-\eta)a_U$.
Since $U\subset\Omega$, the zero extension of each $\mathbf v_U$ lies
in $H_0^1(\Omega)^d$, and its weak gradient is the zero extension of
$\nabla\mathbf v_U$. Using the same notation for these extensions,
set $\mathbf v=\sum_Ua_U\mathbf v_U$.
The overlap bound gives
\[
 \|\nabla\mathbf v\|_{L^2(\Omega)}^2\le M\sum_Ua_U^2,
 \qquad
 \int_\Omega q\,\nabla\cdot\mathbf v\,dx
 \ge(1-\eta)\sum_Ua_U^2.
\]
The definition of the global negative norm implies
\[
 (1-\eta)\sum_Ua_U^2
 \le\|\nabla q\|_{-1,\Omega}\Bigl(M\sum_Ua_U^2\Bigr)^{1/2}.
\]
Letting $\eta\downarrow0$ proves the estimate for every finite
subfamily. Increasing finite subfamilies exhaust the countable case.
\end{proof}

\subsection{Proof of the main theorem}\label{subsec:weighted_necas}

\begin{proof}[Proof of Theorem~\ref{thm:weighted_necas}]
Take the covering and partition of unity from
Lemma~\ref{lem:whitney_covering}. For $Q\in\mathcal W$, let
$m_Q=|U_Q|^{-1}\int_{U_Q}q\,dx$. Lemmas~\ref{lem:cube_necas}
and~\ref{lem:dual_localization} give
\begin{equation}\label{eq:whitney_oscillation}
 \sum_Q\|q-m_Q\|_{L^2(U_Q)}^2
 \le C_d\sum_Q\|\nabla q\|_{-1,U_Q}^2
 \le C_d\|\nabla q\|_{-1,\Omega}^2.
\end{equation}
Define
\begin{equation}\label{eq:whitney_reconstruction}
 \phi:=\sum_{Q\in\mathcal W}m_Q\psi_Q.
\end{equation}
Local finiteness implies $\phi\in C^\infty(\Omega)$.
Since $\sum_Q\psi_Q=1$ and $\sum_Q\nabla\psi_Q=0$, almost everywhere
in $\Omega$ we have
\begin{equation}\label{eq:whitney_cancellation}
 q-\phi=\sum_Q\psi_Q(q-m_Q),\qquad
 \nabla\phi=\sum_Q\nabla\psi_Q(m_Q-q).
\end{equation}
The second identity uses only cancellation of the partition gradients
and does not differentiate $q$. Jensen's inequality gives
\[
 |q-\phi|^2\le\sum_Q\psi_Q|q-m_Q|^2
 \le\sum_Q\chi_{U_Q}|q-m_Q|^2.
\]
At each point at most $N_d$ summands are nonzero. The bound
$\delta|\nabla\psi_Q|\le B_d$ and the Cauchy--Schwarz inequality yield
\[
 |\delta\nabla\phi|^2
 \le N_d B_d^2\sum_Q\chi_{U_Q}|q-m_Q|^2.
\]
Integrating these two inequalities, using Tonelli's theorem and
\eqref{eq:whitney_oscillation}, gives
\begin{equation}\label{eq:whitney_reconstruction_bound}
 \|q-\phi\|_{L^2(\Omega)}^2
 +\|\delta\nabla\phi\|_{L^2(\Omega)}^2
 \le C_d\sum_Q\|q-m_Q\|_{L^2(U_Q)}^2
 \le C_d\|\nabla q\|_{-1,\Omega}^2.
\end{equation}
In particular, $\phi\in L^2(\Omega)$ and
$\delta\nabla\phi\in L^2(\Omega)^d$, so
$\phi\in\mathcal H_\delta(\Omega)$ is an admissible competitor
for the infimum in \eqref{eq:K_delta_def}.
Taking square roots in \eqref{eq:whitney_reconstruction_bound} proves
\eqref{eq:weighted_necas}.
\end{proof}

The theorem requires neither boundary regularity nor a Hardy inequality.
The norm equivalence follows by combining it with the Hardy
reverse estimate below.

\subsection{Norm equivalence under the Hardy inequality}
\label{subsec:norm_equivalence}
\label{subsec:hardy}

The main theorem supplies one side of the equivalence between the
distance-weighted norm and the negative norm of the gradient on every
bounded connected open set. To obtain the reverse estimate, we now
introduce the boundary Hardy inequality, which holds in particular on
bounded Lipschitz domains.

\begin{lemma}[Boundary Hardy inequality on Lipschitz domains]
\label{lem:boundary_hardy}
Let $\Omega\subset\R^d$ be a bounded Lipschitz domain. There exists a
constant $C_H>0$, depending only on $d$ and on the Lipschitz character of
$\Omega$, such that
\begin{equation}\label{eq:hardy}
    \Bigl\|\tfrac{u}{\delta}\Bigr\|_{L^2(\Omega)}
    \;\le\;
    C_H\,\|\nabla u\|_{L^2(\Omega)},
    \qquad
    \forall u\in H_0^1(\Omega).
\end{equation}
\end{lemma}

This is the case $p=2$, $\beta=0$ of the classical weighted Hardy
inequality recalled in \cite[Theorem~1.1]{Lehrbck2014}, extended from
$C_c^\infty(\Omega)$ to $H_0^1(\Omega)$ by density; see also
\cite{grisvard2011elliptic,Egert2015}.

\begin{lemma}[Hardy reverse estimate]\label{lem:hardy_reverse}
Let $\Omega\subset\R^d$ be a bounded connected open set satisfying
\eqref{eq:hardy} with constant $C_H$. Then
\begin{equation}\label{eq:hardy_reverse}
    \|\nabla q\|_{-1,\Omega}
    \le \sqrt{1+C_H^2}\,\|q\|_{K_\delta},
    \qquad q\in L_0^2(\Omega).
\end{equation}
\end{lemma}

\begin{proof}
Let $q\in L_0^2(\Omega)$, fix $\phi\in\mathcal H_\delta(\Omega)$, and
first take $\mathbf v\in C_c^\infty(\Omega)^d$.
Since $\phi\in H^1_{\mathrm{loc}}(\Omega)$, integration by parts gives
\begin{equation*}
    \langle\nabla q,\mathbf v\rangle
    \;=\;
    -\int_\Omega (q-\phi)\,\nabla\!\cdot\!\mathbf v\,dx
    \;+\;
    \int_\Omega \nabla\phi\cdot\mathbf v\,dx.
\end{equation*}
Hence
\begin{equation*}
    |\langle\nabla q,\mathbf v\rangle|
    \;\le\;
    \|q-\phi\|_{L^2(\Omega)}\,\|\nabla\!\cdot\!\mathbf v\|_{L^2(\Omega)}
    \;+\;
    \|\delta\nabla\phi\|_{L^2(\Omega)}\,\Bigl\|\tfrac{\mathbf v}{\delta}\Bigr\|_{L^2(\Omega)}.
\end{equation*}
For $\mathbf v\in H_0^1(\Omega)^d$, the div--curl identity reads
\begin{equation*}
 \|\nabla\mathbf v\|_{L^2(\Omega)}^2
 =\|\nabla\cdot\mathbf v\|_{L^2(\Omega)}^2
  +\frac12\sum_{i,j=1}^d
       \|\partial_jv_i-\partial_iv_j\|_{L^2(\Omega)}^2.
\end{equation*}
In particular,
$\|\nabla\cdot\mathbf v\|_{L^2(\Omega)}
\le\|\nabla\mathbf v\|_{L^2(\Omega)}$.
Using this estimate and the boundary Hardy inequality
\eqref{eq:hardy} componentwise,
\begin{equation*}
    \Bigl\|\tfrac{\mathbf v}{\delta}\Bigr\|_{L^2(\Omega)}
    \;\le\;
    C_H\,\|\nabla\mathbf v\|_{L^2(\Omega)},
\end{equation*}
we obtain
\begin{equation*}
    |\langle\nabla q,\mathbf v\rangle|
    \;\le\;
    \Bigl(\|q-\phi\|_{L^2(\Omega)}
          +C_H\,\|\delta\nabla\phi\|_{L^2(\Omega)}\Bigr)
    \|\nabla\mathbf v\|_{L^2(\Omega)}.
\end{equation*}
By density of $C_c^\infty(\Omega)^d$ in $H_0^1(\Omega)^d$ and continuity
of the gradient pairing, the same bound holds for all
$\mathbf v\in H_0^1(\Omega)^d$.
Dividing by $\|\nabla\mathbf v\|_{L^2(\Omega)}$ and taking the supremum
over $\mathbf v\ne\mathbf0$ gives a bound for $\|\nabla q\|_{-1,\Omega}$.
The Cauchy--Schwarz inequality in $\R^2$, followed by the infimum over
$\phi\in\mathcal H_\delta(\Omega)$, yields
\begin{equation*}
    \|\nabla q\|_{-1,\Omega}
    \;\le\;
    \sqrt{1+C_H^2}\,\|q\|_{K_\delta}.
\end{equation*}
This proves \eqref{eq:hardy_reverse}.
\end{proof}

\begin{corollary}[Distance-weighted norm equivalence]
\label{cor:delta_equivalence}
Under the assumptions of Lemma~\ref{lem:hardy_reverse},
\begin{equation}\label{eq:delta_equivalence}
    c_d\,\|q\|_{K_\delta}
    \;\le\;
    \|\nabla q\|_{-1,\Omega}
    \;\le\;
    \sqrt{1+C_H^2}\,\|q\|_{K_\delta},
    \qquad q\in L_0^2(\Omega),
\end{equation}
where $c_d=C_d^{-1}>0$ depends only on $d$.
In particular, the equivalence holds on bounded connected Lipschitz
domains, with the upper constant controlled by their Hardy constants.
\end{corollary}

\begin{proof}
Combine Theorem~\ref{thm:weighted_necas} with
Lemma~\ref{lem:hardy_reverse}.
\end{proof}

\begin{remark}[Uniform Hardy constants]
\label{rem:uniform_hardy}
The Hardy constant admits geometry-independent bounds in several important
classes. For bounded convex domains in any dimension, one may take
$C_H\le2$ \cite{LaptevSobolev2008}. For bounded simply connected planar
domains, Ancona's result
gives $C_H\le4$ \cite{Ancona1986}; see also Laptev--Sobolev
\cite{LaptevSobolev2008}, who obtain improved bounds for geometrically
restricted classes of simply connected planar domains. The bound $C_H\le4$ is unaffected
by variations in channel width, bends, or constrictions, provided the
planar domain remains simply connected.

For sufficient conditions ensuring uniform Hardy constants in multiply
connected domains and in higher dimensions, see Korte--Lehrb\"ack--Tuominen
\cite{KorteLehrbackTuominen2011}.
\end{remark}

\section{Application to Stokes Preconditioning}\label{sec:preconditioning}

\subsection{The Stokes problem and pressure Schur complement}
\label{subsec:stokes_infsup}
\label{subsec:schur_complement}

We now apply the norm equivalence to preconditioning the steady
incompressible Stokes problem on a bounded connected Lipschitz domain:
\begin{equation}\label{eq:stokes_strong}
    \begin{cases}
        -\mu \Delta \mathbf{u} + \nabla p = f & \text{in } \Omega,          \\
        \nabla \cdot \mathbf{u} = 0           & \text{in } \Omega,          \\
        \mathbf{u} = \mathbf u_D              & \text{on } \partial \Omega.
    \end{cases}
\end{equation}
Here $\mu>0$ is the viscosity, $f$ is the body force, and $\mathbf u_D$
is prescribed velocity data with zero net boundary flux. We set $\mu=1$
in all experiments. After lifting the boundary data, the variational
formulation uses the velocity space $V(\Omega)=H_0^1(\Omega)^d$ with norm
$\|\mathbf v\|_V=\|\nabla\mathbf v\|_{L^2(\Omega)}$ and the pressure
space $Q(\Omega)=L_0^2(\Omega)$. The classical Ne\v{c}as inequality
\eqref{eq:intro_classical_necas} is equivalent to the Stokes inf-sup
condition, with $C_N(\Omega)=\beta(\Omega)^{-1}$ for the optimal
Ne\v{c}as constant $C_N(\Omega)$ and the Stokes inf-sup constant
$\beta(\Omega)$; see \cite{Bramble2003}.

Let $V_h\subset V(\Omega)$ and $Q_{h,0}\subset Q(\Omega)$ be an inf-sup
stable conforming finite-element pair. For unit viscosity, let $A_h$ be the matrix of
$(\nabla\mathbf u_h,\nabla\mathbf v_h)$ and $B_h$ the matrix of
$-(q_h,\nabla\cdot\mathbf v_h)$. After lifting the boundary data, the
discrete Stokes problem has the form
\begin{equation}\label{eq:stokes_block_system}
 \mathcal A_h
 \begin{bmatrix}\mathbf u\\\mathbf p\end{bmatrix}
 =\begin{bmatrix}\mathbf f_h\\\mathbf g_h\end{bmatrix},
 \qquad
 \mathcal A_h=
 \begin{bmatrix}A_h&B_h^\top\\B_h&0\end{bmatrix}.
\end{equation}
Eliminating the velocity gives the pressure Schur complement
$S_h=B_hA_h^{-1}B_h^\top$. Its quadratic form is
\begin{equation}\label{eq:discrete_schur_norm}
 \mathbf q^\top S_h\mathbf q
 =\sup_{\mathbf v_h\in V_h\setminus\{\mathbf0\}}
   \frac{\left(\int_\Omega q_h\,\nabla\cdot\mathbf v_h\,dx\right)^2}
        {\|\nabla\mathbf v_h\|_{L^2(\Omega)}^2},
 \qquad q_h\in Q_{h,0},
\end{equation}
where $\mathbf q$ is the coefficient vector of $q_h$. Thus a block
solver needs an approximation to $S_h^{-1}$ as well as a velocity solver.
Forming $S_h$ is generally expensive because it contains $A_h^{-1}$.

The corresponding continuous Schur complement $S$ has quadratic form
$\langle Sq,q\rangle=\|\nabla q\|_{-1,\Omega}^2$, by
the definition in \eqref{eq:gradient_negative_norm}. A pressure mass inverse uses the
$L^2$ norm to approximate this form. The classical estimates give
\begin{equation*}
 \beta(\Omega)^2\|q\|_{L^2(\Omega)}^2
 \le\langle Sq,q\rangle\le\|q\|_{L^2(\Omega)}^2,
\end{equation*}
so the mass-based bound can deteriorate as $\beta(\Omega)$ decreases.
For the weighted pressure norm, Corollary~\ref{cor:delta_equivalence} gives
\begin{equation}\label{eq:continuous_schur_K}
 c_d^2\|q\|_{K_\delta}^2
 \le\langle Sq,q\rangle
 \le(1+C_H^2)\|q\|_{K_\delta}^2.
\end{equation}
This identifies the distance-weighted norm as a pressure metric whose
equivalence constants are controlled by the dimension and the Hardy
constant. Approximating the Schur complement by the operator associated
with this norm follows the operator-preconditioning framework of
\cite{Mardal2010,loghin2004analysis}.

\subsection{From the weighted norm to a pressure preconditioner}
\label{subsec:weighted_pressure_block}

Assume $Q_{h,0}=Q_h\cap L_0^2(\Omega)$ for a pressure space
$Q_h\subset H^1(\Omega)$ containing constants, and let $\delta_h$
approximate the boundary distance. Since $\Omega$ is bounded,
$Q_h\subset\mathcal H_\delta(\Omega)$.
In a basis of the mean-zero pressure space,
let $M_p$ and $L_{p,\delta_h}$ be the matrices of
\begin{equation*}
 m(q_h,r_h)=\int_\Omega q_hr_h\,dx,
 \qquad
 \ell_{\delta_h}(q_h,r_h)
 =\int_\Omega\delta_h^2\nabla q_h\cdot\nabla r_h\,dx.
\end{equation*}
All pressure matrices and inverses in this section are restricted to
the mean-zero space; the weighted stiffness matrix is assumed positive
definite there. The stiffness form imposes natural Neumann conditions
on the auxiliary pressure problem.

Restricting the infimum in \eqref{eq:K_delta_def} to $Q_{h,0}$ and
using the discrete distance $\delta_h$ gives
\begin{equation}\label{eq:discrete_K_norm}
 \|q_h\|_{K_{\delta_h,h}}^2
 :=\min_{\phi_h\in Q_{h,0}}
 \left(\|q_h-\phi_h\|_{L^2(\Omega)}^2
       +\|\delta_h\nabla\phi_h\|_{L^2(\Omega)}^2\right).
\end{equation}
Restricting $\phi_h$ to have zero mean does not change this minimum:
removing its mean reduces the first term and leaves the second term
unchanged. In coefficient vectors, \eqref{eq:discrete_K_norm} becomes
\begin{equation*}
 \min_{\boldsymbol\varphi}
 \left((\mathbf q-\boldsymbol\varphi)^\top M_p
                    (\mathbf q-\boldsymbol\varphi)
       +\boldsymbol\varphi^\top L_{p,\delta_h}\boldsymbol\varphi\right).
\end{equation*}
The minimizing vector satisfies
$M_p(\mathbf q-\boldsymbol\varphi)
=L_{p,\delta_h}\boldsymbol\varphi=:\mathbf t$. Hence
$\mathbf q=(M_p^{-1}+L_{p,\delta_h}^{-1})\mathbf t$, and the minimum is
\begin{equation}\label{eq:discrete_K_riesz}
 \|q_h\|_{K_{\delta_h,h}}^2=\mathbf q^\top R_{\delta,h}\mathbf q,
 \qquad
 R_{\delta,h}=(M_p^{-1}+L_{p,\delta_h}^{-1})^{-1}.
\end{equation}
We use this norm matrix as the pressure preconditioner. Its inverse is
the sum
\begin{equation}\label{eq:weighted_inverse_block}
 P_\delta:=R_{\delta,h},\qquad
 P_\delta^{-1}=M_p^{-1}+L_{p,\delta_h}^{-1}.
\end{equation}
This additive inverse can be viewed in the space decomposition and
subspace correction framework of Xu~\cite{Xu1992}, with the two
components $q_h-\phi_h$ and $\phi_h$ measured by the mass and weighted
stiffness energies, respectively.
Applying $P_\delta^{-1}$ amounts to adding the results of a pressure mass
solve and a distance-weighted stiffness solve. The weight uses local
geometric information directly, whereas a constant-scale inverse block
$M_p^{-1}+(cL_p)^{-1}$ assigns one coefficient to the whole domain.

The matrix identity \eqref{eq:discrete_K_riesz} is exact. To obtain
uniform spectral bounds for $P_\delta^{-1}S_h$ from
\eqref{eq:continuous_schur_K}, one additionally needs uniform comparison
of the discrete supremum in \eqref{eq:discrete_schur_norm} with the
continuous negative norm, and of \eqref{eq:discrete_K_norm} with
$K_\delta$. If these comparisons hold uniformly, then
\begin{equation}\label{eq:discrete_schur_comparison}
 c\,\mathbf q^\top R_{\delta,h}\mathbf q
 \le\mathbf q^\top S_h\mathbf q
 \le C\,\mathbf q^\top R_{\delta,h}\mathbf q,
 \qquad
 \kappa(P_\delta^{-1}S_h)\le C/c,
\end{equation}
where $\kappa$ denotes the ratio of the largest to smallest eigenvalue.
The constants depend on the comparison bounds as well as the dimension
and the Hardy constant. Thus the continuous result motivates the
discrete block, while its spectral quality and the effect of approximate
auxiliary solves are assessed in Section~\ref{sec:numerical}.

For the full Stokes system, we use an upper triangular block
preconditioner \cite{loghin2004analysis}, with $P$ approximating
$S_h$:
\begin{equation}\label{eq:stokes_block_preconditioner}
 \mathcal P_U=
 \begin{bmatrix}\hat A&B_h^\top\\0&-P\end{bmatrix},
\end{equation}
where $\hat A$ approximates $A_h$.

\section{Numerical Results}\label{sec:numerical}
\raggedbottom
\newcommand{\CampaignCases}{49}
\newcommand{\CampaignSolves}{686}
\newcommand{\CampaignSpectralWins}{49}
\newcommand{\CampaignMaxResidual}{1.9522\times10^{-10}}
\newcommand{\CampaignExactRange}{14--21}
\newcommand{\CampaignAMGRange}{17--25}
\newcommand{\SlenderStartGain}{1.15}
\newcommand{\SlenderEndGain}{1.12}
\newcommand{\SlenderMassConditions}{20.96\text{--}1246.39}
\newcommand{\SlenderMassStart}{1246.39}
\newcommand{\SlenderMassEnd}{20.96}
\newcommand{\SlenderConstantConditions}{8.13\text{--}8.75}
\newcommand{\SlenderConstantStart}{8.75}
\newcommand{\SlenderConstantEnd}{8.13}
\newcommand{\SlenderTunedConditions}{5.62\text{--}5.97}
\newcommand{\SlenderTunedStart}{5.97}
\newcommand{\SlenderTunedEnd}{5.62}
\newcommand{\SlenderWeightedConditions}{5.01\text{--}5.21}
\newcommand{\SlenderWeightedStart}{5.21}
\newcommand{\SlenderWeightedEnd}{5.01}
\newcommand{\SlenderExactRange}{16--17}
\newcommand{\SlenderAmgRange}{20--23}
\newcommand{\TaperStartGain}{1.04}
\newcommand{\TaperEndGain}{3.07}
\newcommand{\TaperMassConditions}{20.96\text{--}80.01}
\newcommand{\TaperMassStart}{20.96}
\newcommand{\TaperMassEnd}{80.01}
\newcommand{\TaperConstantConditions}{8.59\text{--}19.91}
\newcommand{\TaperConstantStart}{9.33}
\newcommand{\TaperConstantEnd}{19.91}
\newcommand{\TaperTunedConditions}{6.20\text{--}18.69}
\newcommand{\TaperTunedStart}{6.20}
\newcommand{\TaperTunedEnd}{18.69}
\newcommand{\TaperWeightedConditions}{5.97\text{--}6.11}
\newcommand{\TaperWeightedStart}{5.97}
\newcommand{\TaperWeightedEnd}{6.09}
\newcommand{\TaperExactRange}{16--18}
\newcommand{\TaperAmgRange}{23--25}
\newcommand{\DensityStartGain}{2.04}
\newcommand{\DensityEndGain}{2.33}
\newcommand{\DensityMassConditions}{40.61\text{--}292.45}
\newcommand{\DensityMassStart}{40.61}
\newcommand{\DensityMassEnd}{292.45}
\newcommand{\DensityConstantConditions}{12.04\text{--}13.11}
\newcommand{\DensityConstantStart}{12.04}
\newcommand{\DensityConstantEnd}{13.11}
\newcommand{\DensityTunedConditions}{9.78\text{--}11.83}
\newcommand{\DensityTunedStart}{9.78}
\newcommand{\DensityTunedEnd}{11.83}
\newcommand{\DensityWeightedConditions}{4.80\text{--}5.08}
\newcommand{\DensityWeightedStart}{4.80}
\newcommand{\DensityWeightedEnd}{5.08}
\newcommand{\DensityExactRange}{17--18}
\newcommand{\DensityAmgRange}{21--23}
\newcommand{\RadiusStartGain}{1.60}
\newcommand{\RadiusEndGain}{10.91}
\newcommand{\RadiusMassConditions}{28.06\text{--}2574.97}
\newcommand{\RadiusMassStart}{28.06}
\newcommand{\RadiusMassEnd}{2574.97}
\newcommand{\RadiusConstantConditions}{10.26\text{--}202.18}
\newcommand{\RadiusConstantStart}{10.26}
\newcommand{\RadiusConstantEnd}{202.18}
\newcommand{\RadiusTunedConditions}{8.81\text{--}59.14}
\newcommand{\RadiusTunedStart}{8.81}
\newcommand{\RadiusTunedEnd}{59.14}
\newcommand{\RadiusWeightedConditions}{5.03\text{--}5.51}
\newcommand{\RadiusWeightedStart}{5.51}
\newcommand{\RadiusWeightedEnd}{5.42}
\newcommand{\RadiusExactRange}{17--18}
\newcommand{\RadiusAmgRange}{22--25}
\newcommand{\GradedStartGain}{2.94}
\newcommand{\GradedEndGain}{5.28}
\newcommand{\GradedMassConditions}{91.15\text{--}263.88}
\newcommand{\GradedMassStart}{263.88}
\newcommand{\GradedMassEnd}{91.15}
\newcommand{\GradedConstantConditions}{27.30\text{--}30.03}
\newcommand{\GradedConstantStart}{27.30}
\newcommand{\GradedConstantEnd}{30.03}
\newcommand{\GradedTunedConditions}{15.15\text{--}29.63}
\newcommand{\GradedTunedStart}{15.15}
\newcommand{\GradedTunedEnd}{29.63}
\newcommand{\GradedWeightedConditions}{5.16\text{--}5.61}
\newcommand{\GradedWeightedStart}{5.16}
\newcommand{\GradedWeightedEnd}{5.61}
\newcommand{\GradedExactRange}{18--18}
\newcommand{\GradedAmgRange}{23--23}
\newcommand{\CompositeGapStartGain}{15.13}
\newcommand{\CompositeGapEndGain}{2.31}
\newcommand{\CompositeGapMassConditions}{175.43\text{--}6952.51}
\newcommand{\CompositeGapMassStart}{6952.51}
\newcommand{\CompositeGapMassEnd}{175.43}
\newcommand{\CompositeGapConstantConditions}{41.83\text{--}87.02}
\newcommand{\CompositeGapConstantStart}{87.02}
\newcommand{\CompositeGapConstantEnd}{43.08}
\newcommand{\CompositeGapTunedConditions}{13.36\text{--}87.02}
\newcommand{\CompositeGapTunedStart}{87.02}
\newcommand{\CompositeGapTunedEnd}{13.36}
\newcommand{\CompositeGapWeightedConditions}{5.71\text{--}5.79}
\newcommand{\CompositeGapWeightedStart}{5.75}
\newcommand{\CompositeGapWeightedEnd}{5.79}
\newcommand{\CompositeGapExactRange}{20--21}
\newcommand{\CompositeGapAmgRange}{24--25}
\newcommand{\CompositeDensityStartGain}{4.98}
\newcommand{\CompositeDensityEndGain}{4.44}
\newcommand{\CompositeDensityMassConditions}{1022.74\text{--}1035.99}
\newcommand{\CompositeDensityMassStart}{1022.74}
\newcommand{\CompositeDensityMassEnd}{1035.99}
\newcommand{\CompositeDensityConstantConditions}{27.10\text{--}67.60}
\newcommand{\CompositeDensityConstantStart}{27.10}
\newcommand{\CompositeDensityConstantEnd}{67.60}
\newcommand{\CompositeDensityTunedConditions}{25.33\text{--}27.10}
\newcommand{\CompositeDensityTunedStart}{27.10}
\newcommand{\CompositeDensityTunedEnd}{26.10}
\newcommand{\CompositeDensityWeightedConditions}{5.44\text{--}5.87}
\newcommand{\CompositeDensityWeightedStart}{5.44}
\newcommand{\CompositeDensityWeightedEnd}{5.87}
\newcommand{\CompositeDensityExactRange}{20--20}
\newcommand{\CompositeDensityAmgRange}{24--24}
\newcommand{\CompositeRadiusStartGain}{4.46}
\newcommand{\CompositeRadiusEndGain}{4.01}
\newcommand{\CompositeRadiusMassConditions}{1013.73\text{--}1043.01}
\newcommand{\CompositeRadiusMassStart}{1043.01}
\newcommand{\CompositeRadiusMassEnd}{1013.73}
\newcommand{\CompositeRadiusConstantConditions}{41.12\text{--}50.59}
\newcommand{\CompositeRadiusConstantStart}{50.59}
\newcommand{\CompositeRadiusConstantEnd}{41.12}
\newcommand{\CompositeRadiusTunedConditions}{24.16\text{--}27.73}
\newcommand{\CompositeRadiusTunedStart}{27.73}
\newcommand{\CompositeRadiusTunedEnd}{24.16}
\newcommand{\CompositeRadiusWeightedConditions}{5.71\text{--}6.21}
\newcommand{\CompositeRadiusWeightedStart}{6.21}
\newcommand{\CompositeRadiusWeightedEnd}{6.03}
\newcommand{\CompositeRadiusExactRange}{20--20}
\newcommand{\CompositeRadiusAmgRange}{24--25}

\newcommand{\CompositeUniformGapGain}{16.75}
\newcommand{\CompositeUniformThroatMinimum}{0.541840}
\newcommand{\CompositeUniformRefinedMinimum}{0.543}

\newcommand{\NumericalSpectralCases}{39}
\newcommand{\NumericalSolverCases}{91}
\newcommand{\NumericalNewCases}{81}
\newcommand{\NumericalMaxDofs}{10,136,708}

\newcommand{\TaperSpectralDofs}{4,916--11,903}
\newcommand{\TaperIterationDofs}{11,006--44,041}
\newcommand{\TaperRefinementDofs}{43,164--173,312}
\newcommand{\TaperMinimumRange}{0.53717--0.54043}
\newcommand{\TaperMinimumAngle}{33.63}
\newcommand{\TaperMaximumEdgeRatio}{1.81}
\newcommand{\TaperRefinedEqualWeighted}{23}
\newcommand{\TaperRefinedEqualConstant}{25}
\newcommand{\TaperRefinedSevereWeighted}{27}
\newcommand{\TaperRefinedSevereConstant}{49}

We begin with the unit square to verify the basic performance of the
distance-weighted pressure block from Section~\ref{sec:preconditioning}
and examine its sensitivity to mesh refinement. Slender channels and
periodic pillar arrays provide benchmarks with a natural global length
scale, for which existing analyses motivate scale-dependent pressure
preconditioners \cite{sande2025robust,xin2026uniform}. Comparisons with
the pressure mass matrix and a constant-scale sum-space block test
whether the boundary-distance weight recovers the effectiveness of these
scale-specific choices. To examine local adaptation, we extend these
settings to tapered channels, arrays with nonuniform pillar counts, and
composite porous--channel domains. These variants test whether the same
weight remains effective when different parts of a domain have competing
geometric scales, without choosing a new global scale for each
configuration. Figures~\ref{fig:num_geometry_parameters}
and~\ref{fig:num_composite_geometry} illustrate the test geometries.

\subsection{Common experimental setup}\label{subsec:num_setup}

\subsubsection{Test geometries}\label{subsec:num_geometries}

The square baseline $(0,1)^2$ provides a mesh-refinement control.
Slender channels test aspect-ratio effects on pressure stability
\cite{Dobrowolski2003} and width-scaled preconditioning
\cite{sande2025robust}. Tapered channels extend the test to variable
local widths, as considered in
\cite[Remark~3 and Example~5]{sande2025robust}.
Both channel families are described by
\begin{equation}\label{eq:num_channel_geometry}
 \Omega_{\varepsilon,\rho}
 =\{(x,y):0<x<1,\ 0<y<\varepsilon[1-(1-\rho^{-1})x]\},
\end{equation}
where $\varepsilon$ is the maximum width and $\rho\ge1$ is the ratio of
maximum to minimum width. The case $\rho=1$ gives the rectangular
slender channel $(0,1)\times(0,\varepsilon)$.

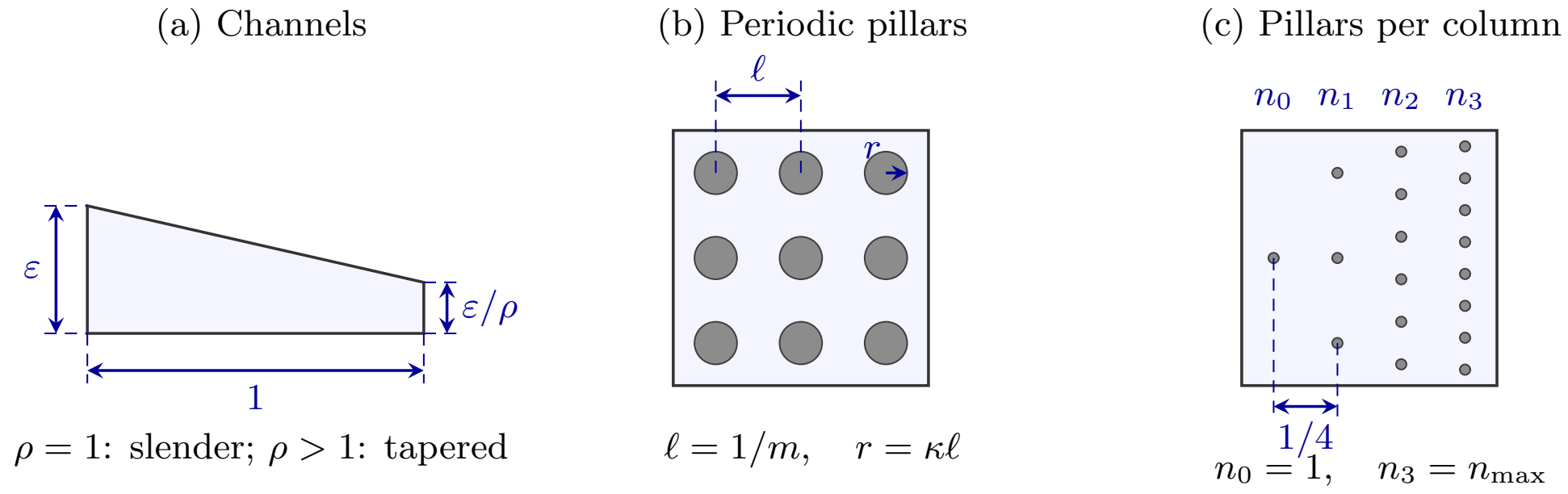
\begin{figure}[!htp]
\centering
\resizebox{\textwidth}{!}{
\begin{tikzpicture}[
    x=1cm,y=1cm,every node/.style={font=\small},
    fluid/.style={draw=black!80,fill=blue!4,line width=0.7pt},
    pillar/.style={draw=black!75,fill=black!45,line width=0.4pt},
    dim/.style={<->,>=stealth,draw=blue!60!black,line width=0.7pt},
    guide/.style={draw=blue!60!black,densely dashed,line width=0.4pt},
    annotation/.style={text=blue!60!black}
]
\begin{scope}
  \node at (2.15,3.2) {(a) Channels};
  \draw[fluid] (0.65,0.55) -- (3.55,0.55) -- (3.55,0.99)
    -- (0.65,1.65) -- cycle;
  \draw[guide] (0.28,0.55) -- (0.65,0.55);
  \draw[guide] (0.28,1.65) -- (0.65,1.65);
  \draw[dim] (0.38,0.55) -- node[left,annotation] {$\varepsilon$} (0.38,1.65);
  \draw[guide] (3.55,0.55) -- (3.85,0.55);
  \draw[guide] (3.55,0.99) -- (3.85,0.99);
  \draw[dim] (3.75,0.55) -- node[right,annotation] {$\varepsilon/\rho$} (3.75,0.99);
  \draw[guide] (0.65,0.55) -- (0.65,0.12);
  \draw[guide] (3.55,0.55) -- (3.55,0.12);
  \draw[dim] (0.65,0.23) -- node[below,annotation] {$1$} (3.55,0.23);
  \node at (2.15,-0.45) {$\rho=1$: slender; $\rho>1$: tapered};
\end{scope}

\begin{scope}[xshift=4.9cm]
  \node at (2.0,3.2) {(b) Periodic pillars};
  \draw[fluid] (0.8,0.1) rectangle (3.0,2.3);
  \foreach \i in {0,...,2} {
    \foreach \j in {0,...,2} {
      \draw[pillar] ({0.8+(\i+0.5)*2.2/3},{0.1+(\j+0.5)*2.2/3})
        circle ({0.25*2.2/3});
    }
  }
  \draw[guide] ({0.8+1.1/3},{0.1+5.5/3}) -- ({0.8+1.1/3},2.73);
  \draw[guide] (1.9,{0.1+5.5/3}) -- (1.9,2.73);
  \draw[dim] ({0.8+1.1/3},2.60) -- node[above,annotation] {$\ell$} (1.9,2.60);
  \draw[dim,->] ({0.8+5.5/3},{0.1+5.5/3}) --
    node[above left,annotation] {$r$} ({0.8+5.5/3+0.25*2.2/3},{0.1+5.5/3});
  \node at (2.0,-0.45) {$\ell=1/m,\quad r=\kappa\ell$};
\end{scope}

\begin{scope}[xshift=9.8cm]
  \node at (2.0,3.2) {(c) Pillars per column};
  \draw[fluid] (0.8,0.1) rectangle (3.0,2.3);
  \foreach \i/\count in {0/1,1/3,2/6,3/8} {
    \foreach \j in {1,...,\count} {
      \draw[pillar] ({0.8+(\i+0.5)*0.55},{0.1+2.2*(\j-0.5)/\count})
        circle (0.045);
    }
    \node[annotation] at ({0.8+(\i+0.5)*0.55},2.57) {$n_{\i}$};
  }
  \draw[guide] (1.075,1.2) -- (1.075,-0.18);
  \draw[guide] (1.625,{0.1+2.2/6}) -- (1.625,-0.18);
  \draw[dim] (1.075,-0.08) -- node[below,annotation] {$1/4$} (1.625,-0.08);
  \node at (2.0,-0.63) {$n_0=1,\quad n_3=n_{\max}$};
\end{scope}
\end{tikzpicture}}
\caption{Construction parameters for channels and pillar arrays.
Grey circles are solid pillars; schematics are not to scale.}
\label{fig:num_geometry_parameters}
\end{figure}

Periodic pillar arrays test the density-dependent stability and
pore-scale pressure norms studied in
\cite{infsup_perforation,xin2026uniform}.
They are obtained by removing $m^2$ circular pillars from
the unit square. Their centres are
$((i+\tfrac12)/m,(j+\tfrac12)/m)$, $0\le i,j<m$, and their radius is
$r=\kappa/m$, where $\kappa$ is the radius relative to the period $1/m$.
The graded-density array tests spatially varying pore scales, a
feature of heterogeneous rock geometries \cite{Mostaghimi2013}.
It also lies in the unit square, with equal-radius
pillars ($r=0.005$) in four columns at $x_i=(i+\tfrac12)/4$,
$i=0,1,2,3$. Column $i$ contains
$n_i=1+\lfloor i(n_{\max}-1)/3+1/2\rfloor$ pillars with vertical
coordinates $y_{ij}=(j+\tfrac12)/n_i$, $j=0,\ldots,n_i-1$.
Thus the pillar count increases from one in the first column to
$n_{\max}$ in the last. Figure~\ref{fig:num_geometry_parameters}
illustrates the channel and pillar-array parameters.

\begin{figure}[!htp]
\centering
\resizebox{0.95\textwidth}{!}{
\begin{tikzpicture}[
    x=1cm,y=1cm,every node/.style={font=\small},
    fluid/.style={draw=black!80,fill=blue!3,line width=0.8pt},
    solid/.style={draw=black!75,fill=black!45,line width=0.5pt},
    dim/.style={<->,>=stealth,draw=blue!60!black,line width=0.7pt},
    guide/.style={draw=blue!60!black,densely dashed,line width=0.45pt},
    annotation/.style={text=blue!60!black}
]
\draw[fluid] (0,0) rectangle (12.9,6);
\fill[blue!9] (0.6,0) rectangle (6.6,6);
\draw[black!45,dashed,line width=0.45pt] (0.6,0) -- (0.6,6.15);
\draw[black!45,dashed,line width=0.45pt] (6.6,0) -- (6.6,6.15);

\foreach \i in {0,...,2} {
  \foreach \j in {0,...,2} {
    \draw[solid] ({6*(0.10+(\i+0.5)/3)},
      {6*((5+7*\j)/24)}) circle (0.5);
  }
}
\draw[solid] (9.0,0) rectangle (10.2,2.4);
\draw[solid] (9.0,3.6) rectangle (10.2,6);
\draw[black!80,line width=0.8pt] (0,0) rectangle (12.9,6);

\node at (7.8,6.47) {Chamber};
\node at (9.6,6.47) {Throat};
\node[right] at (13.1,4.5) {Outlet};
\node at (3.6,3.8) {$m\times m$ pillars};

\draw[guide] (1.6,4.75) -- (1.6,7.13);
\draw[guide] (3.6,4.75) -- (3.6,7.13);
\draw[dim] (1.6,7.0) -- node[above,annotation] {$1/m$} (3.6,7.0);
\node[fill=white,inner sep=1.5pt] at (3.6,6.47) {Porous region};

\draw[dim,->] (1.6,4.75) -- (1.6,5.25);
\fill[blue!60!black] (1.6,4.75) circle (0.8pt);
\node[annotation] at (2.6,5.65) {$r$};
\draw[blue!60!black,line width=0.5pt]
  (2.35,5.5) -- (1.61,5.1);
\draw[dim] (1.6,0.03) -- node[left,annotation] {$w$} (1.6,0.75);
\draw[dim] (1.6,1.75) -- node[left,annotation] {$w$} (1.6,2.5);

\draw[dim] (9.6,2.4) -- node[right,annotation] {$g$} (9.6,3.6);
\draw[black!45,densely dotted,line width=0.45pt] (10.4,3) -- (13.25,3);
\node[right] at (13.25,3) {$y=1/2$};

\draw[->,>=stealth,line width=0.55pt] (0,-0.43) -- (13.3,-0.43) node[right] {$x$};
\foreach \x/\label in {0/0,0.6/0.10,6.6/1.10,9.0/1.50,10.2/1.70,12.9/2.15} {
  \draw[line width=0.5pt] (\x,-0.37) -- (\x,-0.49);
  \node[below] at (\x,-0.49) {$\label$};
}
\draw[->,>=stealth,line width=0.55pt] (-0.45,0) -- (-0.45,6.25) node[above] {$y$};
\foreach \y/\label in {0/0,3/{1/2},6/1} {
  \draw[line width=0.5pt] (-0.51,\y) -- (-0.39,\y);
  \node[left] at (-0.51,\y) {$\label$};
}
\end{tikzpicture}}
\caption{Construction of the composite domain. Grey regions are solid.}
\label{fig:num_composite_geometry}
\end{figure}
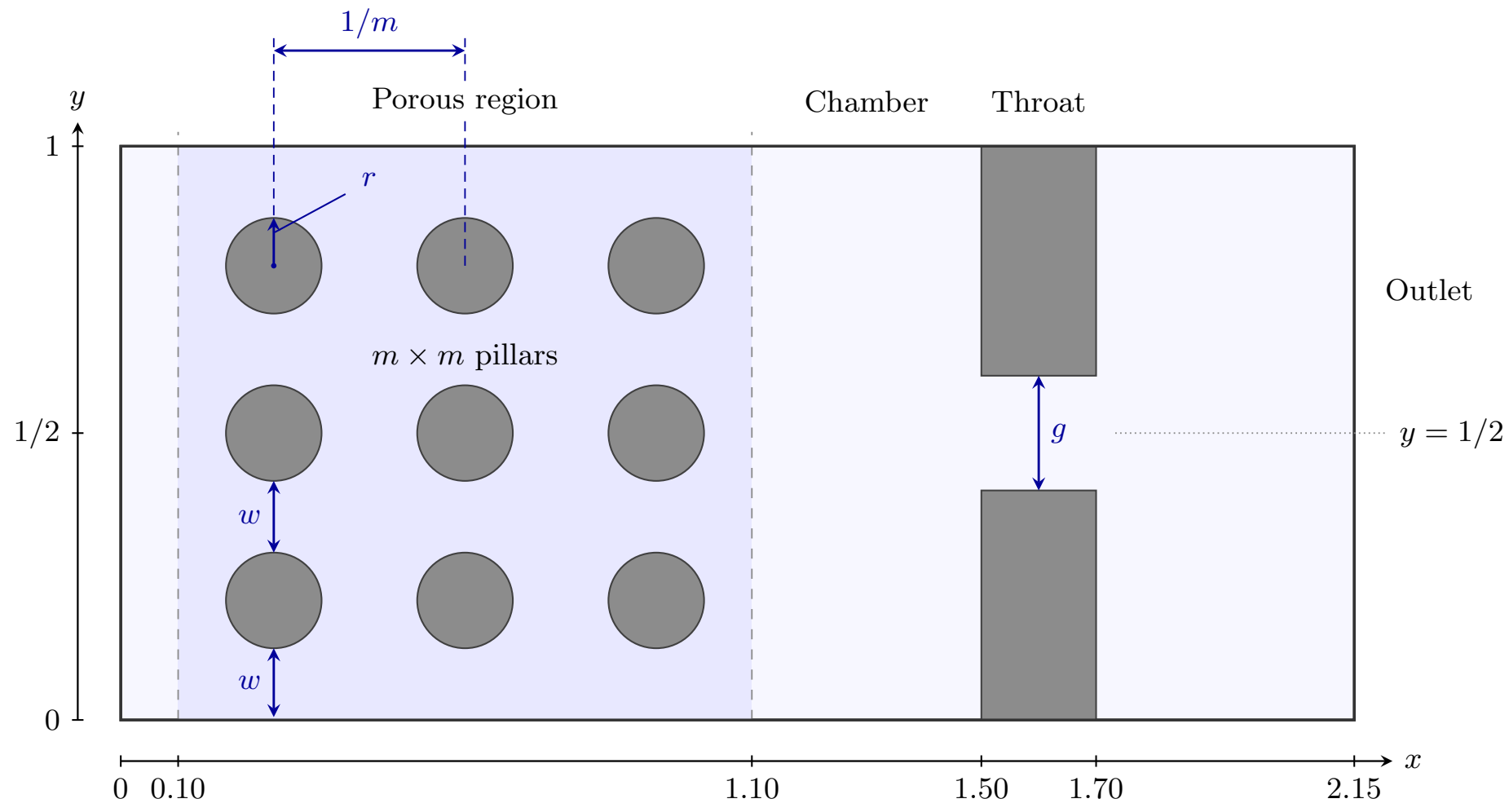

The composite domain is motivated by DLD devices combining pillar
arrays with open collection regions \cite{Loutherback2012}.
We use a simplified geometry with coexisting pore, chamber, and throat
scales. It is formed by removing circular pillars and two
rectangular wall blocks from $(0,2.15)\times(0,1)$. Its square porous region
$(0.10,1.10)\times(0,1)$ contains $m^2$ pillars with centres and radius
\begin{equation*}
 x_{ij}=0.10+(i+\tfrac12)\frac{1}{m},\qquad
 y_{ij}=r+w+j(2r+w),\qquad
 r=\frac{\kappa}{m},
\end{equation*}
where $i,j=0,\ldots,m-1$ and $w=(1-2mr)/(m+1)$, so the vertical
clearances between rows and at the outer walls are equal.
The wall blocks occupy
\begin{equation*}
 [1.50,1.70]\times[0,(1-g)/2],\qquad
 [1.50,1.70]\times[(1+g)/2,1],
\end{equation*}
leaving a throat of width $g$ centred at $y=0.5$. An open chamber
separates the porous region from the throat, followed by an open outlet;
see Figure~\ref{fig:num_composite_geometry} for the construction.
$m^2$ is the pillar count, $\kappa$ is the radius relative to the
horizontal pillar spacing, and $g$ is the throat width.
Parameter ranges and fixed values for each
experiment are given below.

\subsubsection{Discretization and pressure blocks}

All numerical experiments consider the two-dimensional Stokes equations
\eqref{eq:stokes_strong}.
Velocity has prescribed parabolic profiles at the
left and right boundaries, with homogeneous Dirichlet conditions on
all remaining walls, including every obstacle. We use
Taylor--Hood elements with continuous $P_2$ velocity and continuous $P_1$
pressure, restricted to
\begin{equation*}
 Q_{h,0}=\left\{q_h\in Q_h:\int_\Omega q_h\,dx=0\right\}.
\end{equation*}
The boundary distance is approximated by a continuous $P_2$
finite-element function $\delta_h$ on triangular meshes
generated with Gmsh. We use $P_2$ to avoid the identically zero $P_1$
distance interpolant on triangles whose three vertices lie on the boundary.

\begin{figure}[!htp]
\centering
\resizebox{\textwidth}{!}{
\begin{tikzpicture}[
    x=1cm,y=1cm,
    every node/.style={font=\small},
    fluid/.style={draw=black!80,fill=blue!4,line width=0.7pt},
    pillar/.style={draw=black!75,fill=black!45,line width=0.4pt},
    dimension/.style={<->,>=stealth,draw=blue!60!black,line width=0.7pt},
    witness/.style={draw=blue!60!black,densely dashed,line width=0.4pt},
    dimension label/.style={text=blue!60!black,font=\small},
    panel title/.style={font=\small},
    scale label/.style={font=\small}
]

\begin{scope}
  \node[panel title] at (1.95,3.15) {(a) Slender channel};
  \draw[fluid] (0.85,0.7) rectangle (3.85,1.55);
  \draw[witness] (0.32,0.7) -- (0.85,0.7);
  \draw[witness] (0.32,1.55) -- (0.85,1.55);
  \draw[dimension] (0.43,0.7) -- node[left,dimension label] {$W_0$} (0.43,1.55);
  \node[scale label] at (1.95,-0.45) {$W_0=\varepsilon$};
\end{scope}

\begin{scope}[xshift=4.5cm]
  \node[panel title] at (1.95,3.15) {(b) Tapered channel};
  \draw[fluid] (0.85,0.7) -- (3.85,0.7) -- (3.85,0.88)
    -- (0.85,1.55) -- cycle;
  \draw[witness] (0.32,0.7) -- (0.85,0.7);
  \draw[witness] (0.32,1.55) -- (0.85,1.55);
  \draw[dimension] (0.43,0.7) -- node[left,dimension label] {$W_0$} (0.43,1.55);
  \node[scale label] at (1.95,-0.45) {$W_0=\max_x H(x)$};
\end{scope}

\begin{scope}[xshift=9cm]
  \node[panel title] at (1.95,3.15) {(c) Periodic array};
  \draw[fluid] (0.85,0) rectangle (3.05,2.2);
  \foreach \i in {0,...,3} {
    \foreach \j in {0,...,3} {
      \draw[pillar] ({0.85+(\i+0.5)*0.55},{(\j+0.5)*0.55}) circle (0.1375);
    }
  }
  \draw[witness] (1.125,1.925) -- (1.125,2.55);
  \draw[witness] (1.675,1.925) -- (1.675,2.55);
  \draw[dimension] (1.125,2.45) -- node[above,dimension label] {$W_0$} (1.675,2.45);
  \node[scale label] at (1.95,-0.45) {$W_0=\ell=1/m$};
\end{scope}

\begin{scope}[yshift=-4.1cm]
  \node[panel title] at (1.95,3.15) {(d) Graded-density array};
  \draw[fluid] (0.85,0) rectangle (3.05,2.2);
  \foreach \i/\count in {0/1,1/3,2/6,3/8} {
    \foreach \j in {1,...,\count} {
      \draw[pillar] ({0.85+(\i+0.5)*0.55},{2.2*(\j-0.5)/\count}) circle (0.045);
    }
  }
  \draw[witness] (1.125,1.1) -- (1.125,2.55);
  \draw[witness] (1.675,{2.2*5/6}) -- (1.675,2.55);
  \draw[dimension] (1.125,2.45) -- node[above,dimension label] {$W_0$} (1.675,2.45);
  \node[scale label] at (1.95,-0.45) {$W_0=1/4$};
\end{scope}

\begin{scope}[xshift=4.5cm,yshift=-4.1cm]
  \node[panel title] at (1.95,3.15) {(e) Composite domain};
  \draw[fluid] (0.35,0.3) rectangle (3.575,1.8);
  \foreach \i in {0,...,2} {
    \foreach \j in {0,...,2} {
      \draw[pillar] ({0.35+1.5*(0.10+(\i+0.5)/3)},
        {0.3+1.5*(5+7*\j)/24}) circle (0.125);
    }
  }
  \draw[pillar] (2.6,0.3) rectangle (2.9,0.975);
  \draw[pillar] (2.6,1.125) rectangle (2.9,1.8);
  \draw[witness] (0.75,1.4875) -- (0.75,2.55);
  \draw[witness] (1.25,1.4875) -- (1.25,2.55);
  \draw[dimension] (0.75,2.45) -- node[above,dimension label] {$W_0$} (1.25,2.45);
  \node[scale label] at (1.95,-0.45) {$W_0=1/m$};
\end{scope}

\begin{scope}[xshift=9cm,yshift=-4.1cm]
  \node[panel title] at (1.95,3.15) {(f) Square baseline};
  \draw[fluid] (0.85,0) rectangle (3.05,2.2);
  \draw[witness] (0.85,2.2) -- (0.85,2.55);
  \draw[witness] (3.05,2.2) -- (3.05,2.55);
  \draw[dimension] (0.85,2.45) -- node[above,dimension label] {$W_0$} (3.05,2.45);
  \node[scale label] at (1.95,-0.45) {$W_0=1$};
\end{scope}
\end{tikzpicture}}
\caption{Reference lengths $W_0$ (blue arrows) for the constant-scale
pressure block. Grey regions are solid; schematics are not to scale.}
\label{fig:num_reference_scales}
\end{figure}
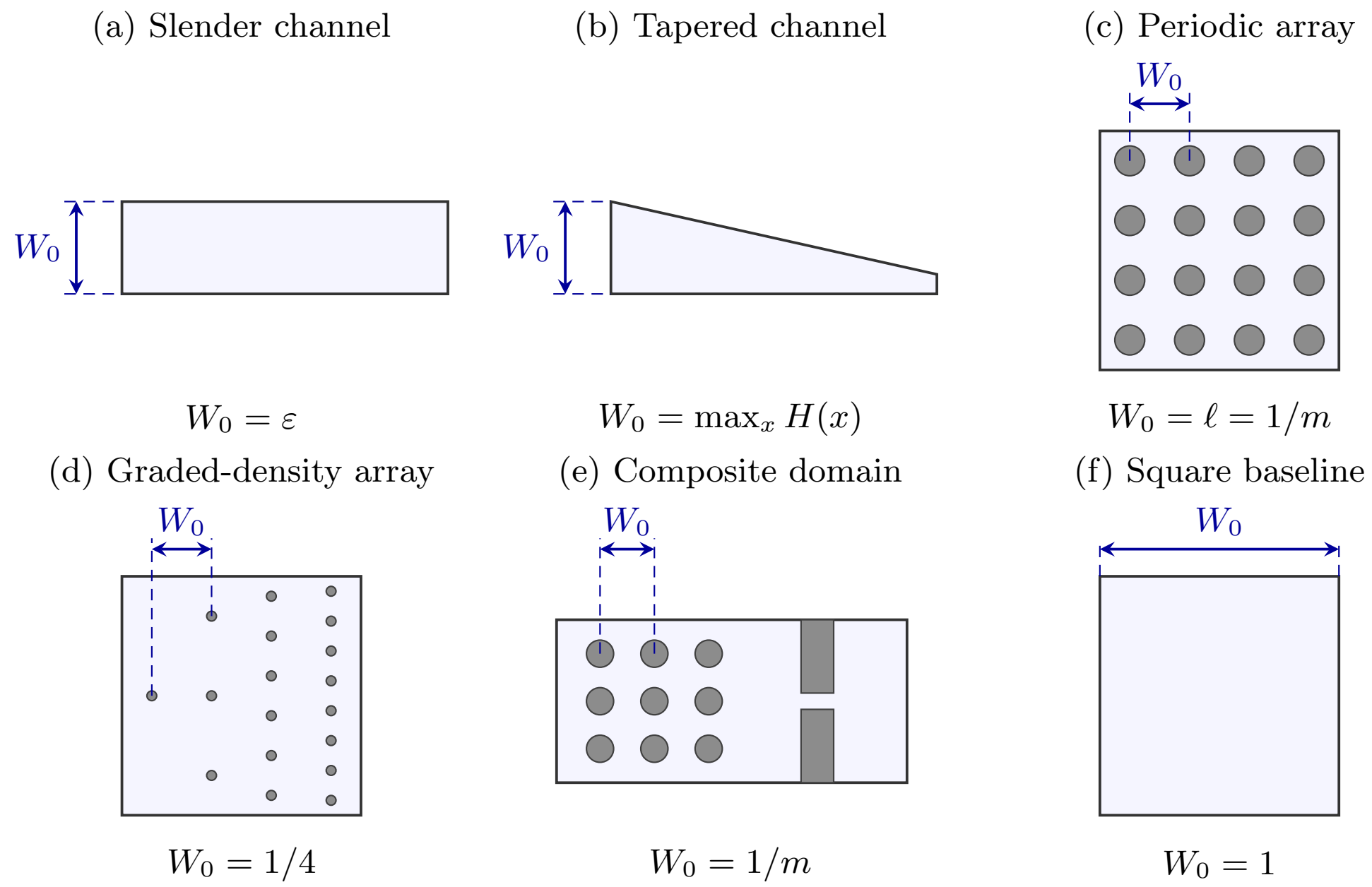

We use the pressure mass matrix $M_p$ and distance-weighted stiffness
matrix $L_{p,\delta_h}$ from Section~\ref{subsec:weighted_pressure_block},
and let $L_p$ denote the unweighted Neumann stiffness matrix.
We compare the pressure preconditioners
\begin{equation}\label{eq:num_pressure_blocks}
 P_M=M_p,\qquad
 P_c=\bigl(M_p^{-1}+(cL_p)^{-1}\bigr)^{-1},\qquad
 P_\delta=\bigl(M_p^{-1}+L_{p,\delta_h}^{-1}\bigr)^{-1}.
\end{equation}
Preconditioning applies $P^{-1}$ for each choice of $P$.
All inverses are interpreted on the mean-zero pressure space, with
natural Neumann conditions for the auxiliary stiffness operators. We
prescribe $c=W_0^2/12$, following the reduced channel model of
Sande et al.~\cite[Section~5]{sande2025robust}. The factor $1/12$ comes
from averaging the parabolic velocity profile between two parallel
no-slip walls.
For channels, $W_0$ is the maximum width; for periodic arrays, it is the
period; for the graded-density array, it is the horizontal distance
between the centre lines of adjacent pillar columns; for the composite domain,
it is the horizontal centre-to-centre pillar spacing $1/m$ in its
porous region; and for the unit-square baseline,
it is one. Figure~\ref{fig:num_reference_scales} illustrates these choices.
The use of this channel coefficient on other geometries
specifies a reproducible baseline, rather than an optimal coefficient.

For a more demanding constant-coefficient comparison, we also compute
$P_{\alpha c}$ for $\alpha\in\{0.1,0.3,1,3,10\}$ on every matrix.
We denote by $P_c^*$ the \emph{best tested constant envelope}.
For each test case, we select the coefficient giving the smallest
condition number among the five candidates. For the FGMRES comparisons
with AMG auxiliary solves, we separately select the coefficient giving
the fewest iterations among the same five candidates.

\subsubsection{Spectral diagnostics}

We examine the eigenproblem $P^{-1}S_h\mathbf q=\lambda\mathbf q$ on
$Q_{h,0}$, where
$S_h=B_hA_h^{-1}B_h^\top$ is the pressure Schur complement. Its equivalent
generalized form is
\begin{equation}\label{eq:num_schur_spectrum}
 S_h\mathbf q=\lambda P\mathbf q.
\end{equation}
For $P=P_M$, the right-hand matrix is $M_p$; for either sum-space block,
\begin{equation*}
 P=(M_p^{-1}+L^{-1})^{-1},\qquad P^{-1}=M_p^{-1}+L^{-1},
 \qquad L=cL_p\ \text{or}\ L_{p,\delta_h}.
\end{equation*}
The smallest and largest eigenvalues of $P^{-1}S_h$ on $Q_{h,0}$ are
denoted by $\lambda_{\min}$ and $\lambda_{\max}$, respectively. Their
Rayleigh quotient characterizations are
\begin{equation*}
 \lambda_{\min}=\min_{\mathbf q\ne\mathbf0}
 \frac{\mathbf q^\top S_h\mathbf q}{\mathbf q^\top P\mathbf q},
 \qquad
 \lambda_{\max}=\max_{\mathbf q\ne\mathbf0}
 \frac{\mathbf q^\top S_h\mathbf q}{\mathbf q^\top P\mathbf q}.
\end{equation*}
Thus $\lambda_{\min}$ and $\lambda_{\max}$ are the optimal lower and
upper constants comparing the Schur complement and preconditioner
quadratic forms. We report their ratio
$\mathcal C(P)=\lambda_{\max}/\lambda_{\min}$ as the spectral condition
number, and list both endpoints for $P_\delta^{-1}S_h$.

For $P=P_\delta$, the denominator is
$\mathbf q^\top P_\delta\mathbf q=\|q_h\|_{K_{\delta_h,h}}^2$ by
\eqref{eq:discrete_K_riesz}. The corresponding continuous comparison
in \eqref{eq:continuous_schur_K} gives
\begin{equation*}
 c_d^2\le
 \frac{\|\nabla q\|_{-1,\Omega}^2}{\|q\|_{K_\delta}^2}
 \le 1+C_H(\Omega)^2,\qquad q\in L_0^2(\Omega)\setminus\{0\}.
\end{equation*}
Thus $\lambda_{\min}$ and $\lambda_{\max}$ correspond to the
geometry-independent lower bound and the Hardy-dependent upper bound,
respectively. We assume that the discrete negative and weighted norms
are uniformly comparable to their continuous counterparts in both mesh
size and geometry, so these bounds carry over up to the comparison constants.

\subsubsection{Prescribed flow and Krylov solver}

There is no body force. On a vertical inlet or outlet of height $H$,
starting at $y=y_0$, we prescribe
\begin{equation}\label{eq:num_parabolic_bc}
 \mathbf u_D(y)=\left(\frac{6F}{H}\eta(1-\eta),0\right)^\top,
 \qquad \eta=\frac{y-y_0}{H}.
\end{equation}
The inlet mean speed is one, so $F=H_{\rm in}$. The same flux $F$ is used
at the outlet to ensure compatibility with incompressibility.
After lifting these boundary data, the velocity correction belongs to
$H_0^1(\Omega)^d$. The lifting changes only the right-hand side, so the
Schur complement uses the homogeneous velocity space introduced in
Section~\ref{subsec:stokes_infsup}.

For each pressure block, the full system is solved with the upper
triangular preconditioner \eqref{eq:stokes_block_preconditioner}.
The exact stage uses sparse direct solves for the auxiliary blocks. The
AMG stage uses one BoomerAMG V-cycle for each inverse application. Both
stages use right-preconditioned FGMRES with a relative PETSc tolerance of
$10^{-10}$; converged runs must also satisfy a measured true relative
residual of $10^{-8}$. All iteration counts refer to outer FGMRES
iterations; the labels Exact and AMG distinguish the auxiliary solvers.
All methods use the same convergence tolerances. NC denotes failure of
FGMRES to converge within 500 iterations.

\subsection{Refinement controls and square baseline}\label{subsec:num_square}

The square $(0,1)^2$ is tested on six mesh levels,
$s=1,2,4,8,16,32$; the slender,
periodic, and composite controls are refined at fixed geometry and $c$
(Table~\ref{tab:num_refinement}). The refinement factor $s$ scales the
target mesh size by $1/s$; $h_{\max}$ is the maximum element diameter,
and $n_q=\dim Q_{h,0}$ is the number of mean-zero pressure degrees of freedom.

\begin{table}[H]
\caption{Refinement at fixed geometry. The slender width is $1/32$;
the periodic parameters are $m=4$, $\kappa=0.25$; the composite parameters
are $m=3$, $\kappa=0.25$, $g=0.1$. The endpoints
$\lambda_{\min}$ and $\lambda_{\max}$ are for $P_\delta^{-1}S_h$;
$\mathcal C(P)$ is the spectral condition number.
Exact and AMG give its FGMRES iteration counts for the
corresponding auxiliary solves.}
\label{tab:num_refinement}
\centering
\resizebox{\textwidth}{!}{%
\begin{tabular}{llrrrrr|rrrrr}
\toprule
Geometry & $s$ & $h_{\max}$ & $n_q$
& $\mathcal C(P_M)$ & $\mathcal C(P_c)$ & $\mathcal C(P_c^*)$ & $\mathcal C(P_\delta)$
& $\boldsymbol{\lambda_{\min}}$ & $\lambda_{\max}$
& Exact & AMG\\
\midrule
Square & 1 & 0.2545 & 43 & 4.46 & 4.43 & 4.35 & 3.36 & 0.534 & 1.798 & 14 & 17 \\
Square & 2 & 0.1214 & 143 & 4.67 & 5.02 & 4.66 & 3.93 & 0.539 & 2.122 & 17 & 18 \\
Square & 4 & 0.0699 & 509 & 4.82 & 5.40 & 4.84 & 4.52 & 0.539 & 2.437 & 17 & 20 \\
Square & 8 & 0.0294 & 1934 & 4.93 & 5.67 & 4.97 & 5.06 & 0.539 & 2.726 & 18 & 21 \\
Square & 16 & 0.0168 & 7546 & 5.01 & 5.87 & 5.07 & 5.53 & 0.540 & 2.983 & 18 & 21 \\
Square & 32 & 0.0081 & 29983 & 5.08 & 6.01 & 5.15 & 5.98 & 0.539 & 3.227 & 18 & 22 \\
Slender & 1 & 0.0094 & 799 & 1246.40 & 7.60 & 5.46 & 4.70 & 0.532 & 2.499 & 16 & 22 \\
Slender & 2 & 0.0051 & 2703 & 1246.39 & 8.75 & 5.97 & 5.21 & 0.539 & 2.806 & 16 & 23 \\
Periodic & 1 & 0.0216 & 4561 & 89.71 & 12.88 & 11.34 & 5.03 & 0.539 & 2.707 & 18 & 22 \\
Periodic & 2 & 0.0105 & 16529 & 90.35 & 12.90 & 12.21 & 5.48 & 0.539 & 2.952 & 18 & 23 \\
Composite & 1 & 0.0347 & 3687 & 1183.17 & 29.49 & 29.49 & 5.23 & 0.543 & 2.838 & 20 & 23 \\
Composite & 2 & 0.0175 & 13952 & 1192.33 & 29.45 & 29.45 & 5.65 & 0.543 & 3.068 & 20 & 24 \\
\bottomrule

\end{tabular}}
\end{table}

The smallest eigenvalues of $P_\delta^{-1}S_h$ in
Table~\ref{tab:num_refinement} remain near $0.539$ across the different
geometries and mesh levels, consistent with the geometry-independent
lower bound in \eqref{eq:continuous_schur_K}. The increase in
$\lambda_{\max}$ under refinement is consistent with richer
finite-element spaces resolving additional pressure modes. On the
convex domains, all computed upper endpoints remain below $5$, the
continuous upper bound obtained from $C_H\le2$. Across all four domain
types, the weighted condition numbers remain modest, ranging from
$3.36$ to $5.98$, and the FGMRES iteration counts remain stable under
refinement with both exact and AMG auxiliary solves.

\subsection{Slenderness and spatially varying width}\label{subsec:num_slender}

For slender channels ($\rho=1$), we vary $\varepsilon$ from $1/4$ to
$1/128$, with $s=2$ and $c=\varepsilon^2/12$. In the spectral subset,
narrowing the channel from $1/4$ to $1/32$ raises $\mathcal C(P_M)$ from
$20.96$ to $1246.39$, while $\mathcal C(P_c)$ remains between $8.13$ and
$8.75$ and $\mathcal C(P_\delta)$ between $5.01$ and $5.21$
(Table~\ref{tab:num_channel_spectra}). Both sum-space blocks therefore
compensate for uniform slenderness.

We then fix $\varepsilon=1/4$ and introduce tapering, with width contrast
$\rho$ up to 8 for the spectral tests ($s=8$) and 48 for the FGMRES tests
($s=12$). As $\rho$ increases from 1 to 8, $\mathcal C(P_c)$ rises from
$9.33$ to $19.91$, and even the tuned $\mathcal C(P_c^*)$ increases from
$6.20$ to $18.69$. In comparison, $\mathcal C(P_\delta)$ stays between
$5.97$ and $6.10$ (Table~\ref{tab:num_channel_spectra}).

A single width provides an effective global scale for a uniform channel;
the distance weight also follows the changing local width in a tapered
channel. Across both families, $\lambda_{\min}(P_\delta^{-1}S_h)$ remains
near $0.539$, consistent with the geometry-independent lower bound.
The FGMRES results in Table~\ref{tab:num_channel_amg} show the same
robustness to both uniform narrowing and spatial width variation. The
mass block deteriorates under narrowing, while the constant-scale blocks
lose robustness as the width contrast increases.

\begin{table}[!htp]
\caption{Channel spectral results. Slender channels use $s=2$; tapered
channels use locally refined physical meshes at $s=8$
and maximum width $1/4$. Taper contrasts are shown at $\rho=1,3,5,8$. The condition numbers compare all
four blocks; the last two columns give the endpoints of $P_\delta^{-1}S_h$.}
\label{tab:num_channel_spectra}
\centering
\resizebox{\textwidth}{!}{%
\begin{tabular}{llrrrr|rrr}
\toprule
Geometry & Parameter & $n_q$ & $\mathcal C(P_M)$ & $\mathcal C(P_c)$
& $\mathcal C(P_c^*)$ & $\mathcal C(P_\delta)$ & $\boldsymbol{\lambda_{\min}}$ & $\lambda_{\max}$\\
\midrule
Slender & $\varepsilon=1/32$ & 2703 & 1246.39 & 8.75 & 5.97 & 5.21 & 0.539 & 2.806 \\
Slender & $\varepsilon=1/16$ & 1417 & 312.68 & 8.71 & 5.96 & 5.21 & 0.539 & 2.806 \\
Slender & $\varepsilon=1/8$ & 682 & 79.25 & 8.59 & 5.88 & 5.15 & 0.539 & 2.779 \\
Slender & $\varepsilon=1/4$ & 361 & 20.96 & 8.13 & 5.62 & 5.01 & 0.539 & 2.697 \\
\midrule
Taper & $\rho=1$ & 4916 & 20.96 & 9.33 & 6.20 & 5.97 & 0.539 & 3.217 \\
Taper & $\rho=3$ & 8252 & 49.31 & 8.59 & 8.59 & 6.09 & 0.539 & 3.278 \\
Taper & $\rho=5$ & 10125 & 65.30 & 12.16 & 12.16 & 6.10 & 0.539 & 3.284 \\
Taper & $\rho=8$ & 11903 & 80.01 & 19.91 & 18.69 & 6.09 & 0.540 & 3.290 \\
\bottomrule

\end{tabular}}
\end{table}

\begin{table}[!htp]
\caption{FGMRES iteration counts for channel problems with one-cycle AMG
auxiliary solves. Slender channels use $s=2$; tapered channels use physical
meshes at $s=12$. Each column is a parameter value; NC denotes failure to
converge within 500 iterations.}
\label{tab:num_channel_amg}
\centering
\small
\setlength{\tabcolsep}{2pt}
\renewcommand{\arraystretch}{1.08}
\begin{tabularx}{\textwidth}{@{}>{\raggedright\arraybackslash}p{0.15\textwidth}*{11}{>{\centering\arraybackslash}X}@{}}
\toprule
Slender $\varepsilon$ & $1/4$ & $1/6$ & $1/8$ & $1/12$ & $1/16$ & $1/24$ & $1/32$ & $1/48$ & $1/64$ & $1/96$ & $1/128$ \\
\midrule
$P_M$ & 29 & 33 & 36 & 45 & 54 & 97 & 155 & 198 & 262 & NC & NC \\
$P_c$ & 24 & 24 & 23 & 26 & 26 & 25 & 25 & 26 & 26 & 26 & 26 \\
$P_c^*$ & 24 & 24 & 23 & 25 & 25 & 25 & 25 & 25 & 25 & 26 & 26 \\
$P_\delta$ & 20 & 21 & 21 & 22 & 23 & 23 & 23 & 23 & 23 & 24 & 24 \\
\bottomrule
\end{tabularx}
\par\medskip
\begin{tabularx}{\textwidth}{@{}>{\raggedright\arraybackslash}p{0.15\textwidth}*{11}{>{\centering\arraybackslash}X}@{}}
\toprule
Taper $\rho$ & $1$ & $2$ & $3$ & $4$ & $6$ & $8$ & $12$ & $16$ & $24$ & $32$ & $48$ \\
\midrule
$P_M$ & 28 & 36 & 38 & 40 & 43 & 45 & 48 & 50 & 53 & 56 & 58 \\
$P_c$ & 25 & 28 & 30 & 31 & 32 & 34 & 36 & 39 & 41 & 44 & 47 \\
$P_c^*$ & 25 & 28 & 30 & 31 & 31 & 32 & 35 & 36 & 39 & 42 & 44 \\
$P_\delta$ & 23 & 24 & 25 & 25 & 25 & 25 & 26 & 25 & 25 & 26 & 25 \\
\bottomrule
\end{tabularx}

\end{table}

\subsection{Obstacle density and obstacle size}\label{subsec:num_periodic}

For the periodic arrays defined in Section~\ref{subsec:num_geometries},
the period is $\ell=1/m$, the radius is $r=\kappa\ell$, and we set
$c=\ell^2/12$. We separate changes in the overall pore scale from
changes in the contrast between local scales. The \emph{density sweep}
fixes $\kappa=0.25$
and uses
\begin{equation*}
 m\in\{2,3,4,6,8,10,12,16,24,32,48,64\}.
\end{equation*}
The obstacle count increases from 4 to 4096, with pillars and passages
shrinking together. The
\emph{size sweep} fixes $m=4$ and varies
$\kappa$ over twelve values between 0.1 and 0.495. As $\kappa$ approaches
$1/2$, the passages between neighbouring pillars become arbitrarily
narrow, while the wider pore regions retain a scale comparable to
$\ell$. This increases the contrast between local boundary-distance
scales within the same domain, which a single period-based coefficient
cannot track.

Table~\ref{tab:num_obstacle_spectra} shows that the smallest eigenvalues
of $P_\delta^{-1}S_h$ remain between $0.538$ and $0.539$ in both sweeps.
Both sum-space blocks remain stable as density increases. As the radius
increases, however, the weighted condition number changes little while
the tuned constant envelope deteriorates. The FGMRES results in
Table~\ref{tab:num_obstacle_amg} show the same distinction: at
$\kappa=0.495$, the weighted block retains stable convergence, while
the mass block and all five constant candidates fail to converge within
500 iterations.

\subsubsection{Increasing pillar counts across columns}\label{subsec:num_graded_density}

For the graded-density array defined in Section~\ref{subsec:num_geometries},
we increase $n_{\max}$ from 1 to 99 while keeping the pillar radius and
column positions fixed. This creates increasingly narrow gaps in the
denser columns. The weighted block remains robust throughout the sweep,
whereas the mass block and all five constant
candidates fail to converge within 500 iterations at $n_{\max}=99$
(Table~\ref{tab:num_obstacle_amg}). These results illustrate the benefit
of tracking local scales when the pillar density varies within one domain.

\begin{table}[!htp]
\caption{Obstacle spectral results. Density varies $m$ at $\kappa=0.25$;
size varies $\kappa$ at $m=4$.
The last two columns are the endpoints of $P_\delta^{-1}S_h$.}
\label{tab:num_obstacle_spectra}
\centering
\resizebox{\textwidth}{!}{%
\begin{tabular}{llrrrr|rrr}
\toprule
Sweep & Parameter & $n_q$ & $\mathcal C(P_M)$ & $\mathcal C(P_c)$
& $\mathcal C(P_c^*)$ & $\mathcal C(P_\delta)$ & $\boldsymbol{\lambda_{\min}}$ & $\lambda_{\max}$\\
\midrule
Density & $m=2$ & 1173 & 40.61 & 12.04 & 9.78 & 4.80 & 0.539 & 2.589 \\
Density & $m=4$ & 4561 & 89.71 & 12.88 & 11.34 & 5.03 & 0.539 & 2.707 \\
Density & $m=8$ & 18112 & 292.45 & 13.11 & 11.83 & 5.08 & 0.539 & 2.737 \\
\midrule
Size & $\kappa=0.1$ & 15394 & 28.06 & 10.26 & 8.81 & 5.51 & 0.538 & 2.965 \\
Size & $\kappa=0.25$ & 4561 & 89.71 & 12.88 & 11.34 & 5.03 & 0.539 & 2.707 \\
Size & $\kappa=0.4$ & 16308 & 606.01 & 54.75 & 21.74 & 5.26 & 0.539 & 2.836 \\
Size & $\kappa=0.45$ & 46561 & 2574.97 & 202.18 & 59.14 & 5.42 & 0.539 & 2.921 \\
\bottomrule

\end{tabular}}
\end{table}

\begin{table}[!htp]
\caption{FGMRES iteration counts for obstacle arrays with one-cycle AMG
auxiliary solves. Each column is a parameter value; NC denotes failure to
converge within 500 iterations. The density and size sweeps respectively
fix $\kappa=0.25$ and $m=4$; the graded-density sweep varies $n_{\max}$
with four columns and equal radii.}
\label{tab:num_obstacle_amg}
\centering
\small
\setlength{\tabcolsep}{2pt}
\begin{tabularx}{\textwidth}{@{}>{\raggedright\arraybackslash}p{0.15\textwidth}*{12}{>{\centering\arraybackslash}X}@{}}
\toprule
Density $m$ & $2$ & $3$ & $4$ & $6$ & $8$ & $10$ & $12$ & $16$ & $24$ & $32$ & $48$ & $64$ \\
\midrule
$P_M$ & 53 & 64 & 73 & 121 & 142 & 153 & 170 & 197 & 379 & NC & NC & NC \\
$P_c$ & 35 & 34 & 37 & 35 & 35 & 35 & 35 & 37 & 36 & 36 & 36 & 35 \\
$P_c^*$ & 27 & 29 & 31 & 32 & 32 & 31 & 31 & 32 & 31 & 32 & 31 & 31 \\
$P_\delta$ & 21 & 22 & 22 & 22 & 23 & 23 & 23 & 23 & 24 & 24 & 24 & 24 \\
\bottomrule
\end{tabularx}
\par\medskip
\begin{tabularx}{\textwidth}{@{}>{\raggedright\arraybackslash}p{0.15\textwidth}*{12}{>{\centering\arraybackslash}X}@{}}
\toprule
Size $\kappa$ & $0.1$ & $0.2$ & $0.25$ & $0.3$ & $0.35$ & $0.4$ & $0.425$ & $0.45$ & $0.46$ & $0.475$ & $0.49$ & $0.495$ \\
\midrule
$P_M$ & 42 & 61 & 73 & 97 & 140 & 182 & 268 & 428 & NC & NC & NC & NC \\
$P_c$ & 28 & 31 & 37 & 43 & 54 & 76 & 111 & 152 & 220 & 270 & NC & NC \\
$P_c^*$ & 28 & 31 & 31 & 33 & 37 & 45 & 57 & 76 & 89 & 132 & 283 & NC \\
$P_\delta$ & 22 & 22 & 22 & 23 & 23 & 23 & 24 & 25 & 26 & 27 & 25 & 25 \\
\bottomrule
\end{tabularx}
\par\medskip
\begin{tabularx}{\textwidth}{@{}>{\raggedright\arraybackslash}p{0.15\textwidth}*{12}{>{\centering\arraybackslash}X}@{}}
\toprule
Graded $n_{\max}$ & $1$ & $2$ & $4$ & $8$ & $16$ & $32$ & $64$ & $80$ & $90$ & $96$ & $98$ & $99$ \\
\midrule
$P_M$ & 30 & 33 & 35 & 40 & 48 & 64 & 83 & 166 & 234 & NC & NC & NC \\
$P_c$ & 28 & 28 & 29 & 30 & 35 & 45 & 63 & 75 & 115 & 202 & 346 & NC \\
$P_c^*$ & 25 & 26 & 28 & 30 & 35 & 45 & 63 & 75 & 115 & 202 & 291 & NC \\
$P_\delta$ & 22 & 23 & 24 & 23 & 24 & 24 & 27 & 27 & 28 & 28 & 29 & 30 \\
\bottomrule
\end{tabularx}

\end{table}

\subsection{Composite porous--channel geometry}\label{subsec:num_composite}

Numerical simulations of microfluidic devices, such as deterministic
lateral displacement (DLD) chips, often involve flow passages of different
widths within the same geometry. We illustrate this feature with a
simplified composite domain defined in Section~\ref{subsec:num_geometries},
combining a porous array, an open chamber, and a narrow passage.
The prescribed constant uses the porous horizontal pitch,
$c=1/(12m^2)$.

\begin{figure}[!htp]
\centering
\includegraphics[width=\textwidth]{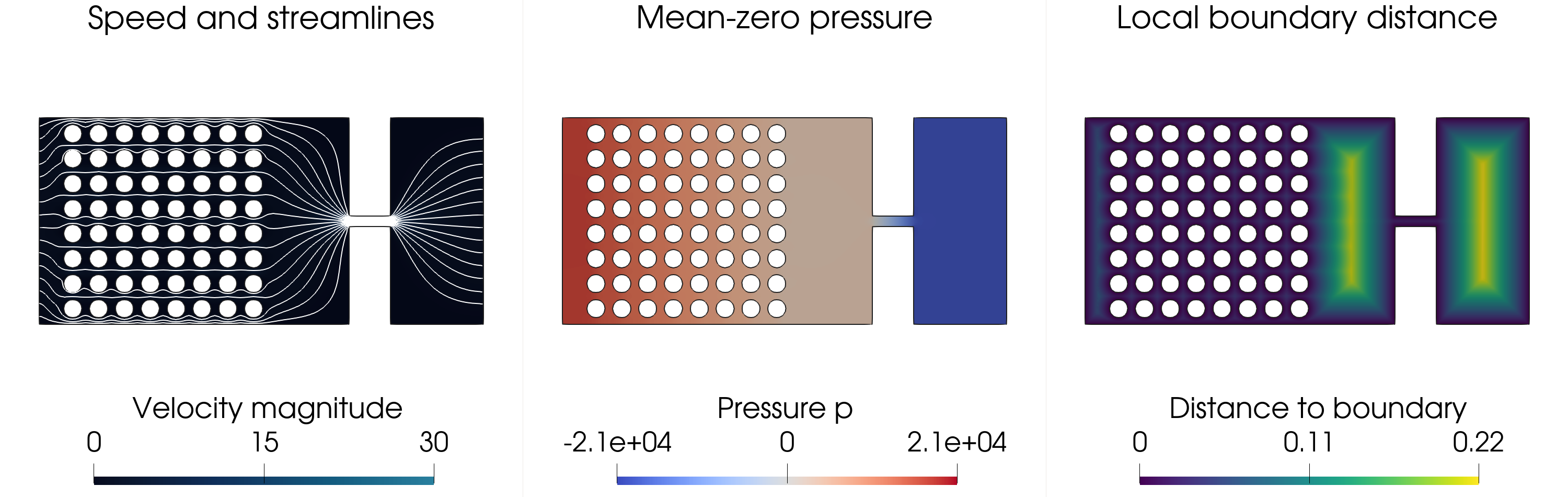}
\caption{Computed dense composite flow at $g=0.05$, $m=8$, $\kappa=0.35$
(64 pillars): velocity magnitude and streamlines (left),
mean-zero pressure (centre), and boundary distance $\delta_h$ (right).}
\label{fig:num_composite_fields}
\end{figure}

Three independent sweeps use twelve values each: (i) $g$ from 0.2 to 0.005 at
$m=3$, $\kappa=0.25$; (ii) $m\in\{2,3,4,5,6,8,10,12,16,24,32,64\}$ at $g=0.1$,
$\kappa=0.25$; and (iii) $\kappa$ from 0.1 to 0.495 at
$m=3$, $g=0.1$. This design tests
whether a local pressure weight can accommodate porous, open, and
narrow-channel scales within the same domain.
Figure~\ref{fig:num_composite_fields} illustrates the computed flow and
local boundary distance on a denser composite with $m=8$, $\kappa=0.35$,
and $g=0.05$.

\begin{table}[!htp]
\caption{FGMRES iteration counts for composite domains with one-cycle AMG
auxiliary solves. Each column is a parameter value; NC denotes failure to
converge within 500 iterations. Fixed parameters are as in
Table~\ref{tab:num_composite_spectra}.}
\label{tab:num_composite_amg}
\centering
\small
\setlength{\tabcolsep}{2pt}
\begin{tabularx}{\textwidth}{@{}>{\raggedright\arraybackslash}p{0.15\textwidth}*{12}{>{\centering\arraybackslash}X}@{}}
\toprule
Throat $g$ & $0.2$ & $0.15$ & $0.1$ & $0.075$ & $0.05$ & $0.04$ & $0.03$ & $0.025$ & $0.02$ & $0.015$ & $0.01$ & $0.005$ \\
\midrule
$P_M$ & 68 & 68 & 72 & 90 & 81 & 113 & 166 & 157 & 207 & 207 & 439 & NC \\
$P_c$ & 33 & 34 & 37 & 39 & 44 & 46 & 52 & 55 & 60 & 70 & 94 & 395 \\
$P_c^*$ & 33 & 34 & 37 & 39 & 44 & 46 & 52 & 55 & 60 & 70 & 94 & 395 \\
$P_\delta$ & 22 & 22 & 23 & 23 & 23 & 24 & 25 & 25 & 25 & 25 & 26 & 27 \\
\bottomrule
\end{tabularx}
\par\medskip
\begin{tabularx}{\textwidth}{@{}>{\raggedright\arraybackslash}p{0.15\textwidth}*{12}{>{\centering\arraybackslash}X}@{}}
\toprule
Density $m$ & $2$ & $3$ & $4$ & $5$ & $6$ & $8$ & $10$ & $12$ & $16$ & $24$ & $32$ & $64$ \\
\midrule
$P_M$ & 60 & 72 & 114 & 168 & 185 & 235 & 298 & 347 & 450 & NC & NC & NC \\
$P_c$ & 36 & 37 & 39 & 41 & 45 & 53 & 60 & 66 & 80 & 117 & 154 & 294 \\
$P_c^*$ & 33 & 37 & 39 & 41 & 45 & 53 & 60 & 66 & 79 & 109 & 145 & 269 \\
$P_\delta$ & 22 & 23 & 23 & 24 & 24 & 24 & 24 & 24 & 25 & 25 & 27 & 29 \\
\bottomrule
\end{tabularx}
\par\medskip
\begin{tabularx}{\textwidth}{@{}>{\raggedright\arraybackslash}p{0.15\textwidth}*{12}{>{\centering\arraybackslash}X}@{}}
\toprule
Size $\kappa$ & $0.1$ & $0.15$ & $0.2$ & $0.25$ & $0.3$ & $0.35$ & $0.4$ & $0.45$ & $0.46$ & $0.475$ & $0.49$ & $0.495$ \\
\midrule
$P_M$ & 56 & 58 & 64 & 72 & 112 & 112 & 238 & NC & NC & NC & NC & NC \\
$P_c$ & 34 & 34 & 35 & 37 & 40 & 50 & 66 & 138 & 140 & 189 & 463 & NC \\
$P_c^*$ & 34 & 34 & 35 & 37 & 40 & 44 & 55 & 87 & 101 & 163 & 376 & NC \\
$P_\delta$ & 24 & 24 & 23 & 23 & 22 & 23 & 23 & 24 & 24 & 24 & 26 & 26 \\
\bottomrule
\end{tabularx}

\end{table}

The composite gives the clearest separation between local and global
scaling because pore gaps, wide chambers, and a narrow throat coexist.
The weighted block remains robust under changes in throat width,
obstacle density, and obstacle radius, while the mass and constant-scale
blocks deteriorate markedly (Table~\ref{tab:num_composite_amg}). At
$\kappa=0.495$, the mass block and all five constant candidates fail to
converge within 500 iterations. The spectral subset also
shows a factor of $\CompositeUniformGapGain$ improvement over the
five-candidate envelope at $g=0.05$.

The refinement control in Table~\ref{tab:num_refinement} supports the
persistence of this robustness under mesh refinement.
Table~\ref{tab:num_composite_spectra} records the corresponding spectral
subset. The weighted lower endpoints again cluster near $0.54$, ranging
from $0.541$ to $0.547$,
while the upper endpoints remain bounded over these sweeps.

\begin{table}[!htp]
\caption{Composite spectral results. The throat sweep fixes $m=3$,
$\kappa=0.25$; the density sweep fixes $g=0.1$, $\kappa=0.25$; the size
sweep fixes $m=3$, $g=0.1$. The last two columns are the endpoints of
$P_\delta^{-1}S_h$. Larger iteration-only cases appear in
Table~\ref{tab:num_composite_amg}.}
\label{tab:num_composite_spectra}
\centering
\resizebox{\textwidth}{!}{%
\begin{tabular}{llrrrr|rrr}
\toprule
Sweep & Parameter & $n_q$ & $\mathcal C(P_M)$ & $\mathcal C(P_c)$
& $\mathcal C(P_c^*)$ & $\mathcal C(P_\delta)$ & $\boldsymbol{\lambda_{\min}}$ & $\lambda_{\max}$\\
\midrule
Throat & $g=0.05$ & 13820 & 7767.48 & 112.95 & 95.24 & 5.69 & 0.544 & 3.093 \\
Throat & $g=0.1$ & 3687 & 1183.17 & 29.49 & 29.49 & 5.23 & 0.543 & 2.838 \\
Throat & $g=0.2$ & 2831 & 283.90 & 15.81 & 15.81 & 5.06 & 0.542 & 2.743 \\
\midrule
Density & $m=2$ & 3654 & 1156.56 & 39.70 & 30.80 & 5.23 & 0.543 & 2.838 \\
Density & $m=3$ & 3687 & 1183.17 & 29.49 & 29.49 & 5.23 & 0.543 & 2.838 \\
Density & $m=4$ & 4949 & 1223.57 & 25.62 & 25.62 & 5.35 & 0.541 & 2.895 \\
\midrule
Size & $\kappa=0.1$ & 17356 & 1185.71 & 32.79 & 32.79 & 5.69 & 0.544 & 3.093 \\
Size & $\kappa=0.25$ & 3687 & 1183.17 & 29.49 & 29.49 & 5.23 & 0.543 & 2.838 \\
Size & $\kappa=0.4$ & 11845 & 2105.64 & 70.71 & 62.53 & 5.63 & 0.547 & 3.083 \\
\bottomrule

\end{tabular}}
\end{table}

\flushbottom
\section{Conclusion}\label{sec:conclusion}

The main result is a distance-weighted Ne\v{c}as-type inequality on every
bounded connected open set:
$\|q\|_{K_\delta}\le C_d\,\|\nabla q\|_{-1,\Omega}$, with $C_d$ depending
only on the dimension. Local Ne\v{c}as estimates on enlarged Whitney
cubes and negative-norm localization control the local oscillations.
A smooth partition of unity combines the cube averages into a competitor
in the natural weighted energy space, with both terms in the
$K$-functional controlled by those oscillations.

Equivalence with the negative norm of the gradient follows as a
corollary when the boundary Hardy inequality holds. The reverse
estimate has constant
$\sqrt{1+C_H(\Omega)^2}$, so the Hardy constant is the only
domain-dependent quantity in the resulting equivalence. The corollary
is uniform on families with uniformly bounded Hardy constants; for all
bounded simply connected planar domains, the reverse constant can be
chosen as $\sqrt{17}$. It motivates a pressure preconditioner that adapts
to the local boundary distance, as examined in Section~\ref{sec:numerical}.

The central numerical finding is that $\lambda_{\min}(P_\delta^{-1}S_h)$
stays close to $0.54$ across the reported geometries and mesh
refinements. This common empirical lower spectral level is consistent
with the geometry-independent lower bound in the continuous estimate.
Across the geometric sweeps, the distance-weighted block consistently
requires about 20--30 FGMRES iterations with AMG auxiliary solves.
The mass and constant-scale blocks lose robustness as the geometry
becomes more extreme; in the near-contact array and composite tests,
the mass block and all five constant candidates fail to converge within
500 iterations. Extending these controls to graded DLD chips and
geological porous media, with scalable weighted auxiliary solves, is a
next step.

\bibliographystyle{siamplain}
\bibliography{reference}

\end{document}